\documentclass[a4paper,12pt,dvips]{article} 
\usepackage{amssymb}
\usepackage{amscd}
\usepackage{amsmath} 
\usepackage{graphicx}
\usepackage{latexsym} 
\usepackage{enumerate}

\usepackage{theorem}
\theoremstyle{plain}
\theorembodyfont{\itshape}
\newtheorem{theorem}{Theorem}[section]

\newtheorem{lemma}[theorem]{Lemma}
 
\theorembodyfont{\rmfamily}

\def\R{\mathbb{R}}
\def\Sol{\mathcal{S}} 
\def\H{\mathcal{H}}
\def\rA{\mathrm{A}}
\def\rB{\mathrm{B}}
\def\rC{\mathrm{C}}
\def\rD{\mathrm{D}}
\def\rO{\mathrm{O}}
\def\rP{\mathrm{P}}
\def\rQ{\mathrm{Q}}
\def\rR{\mathrm{R}}
\def\ra{\mathrm{a}}
\def\rb{\mathrm{b}}
\def\rc{\mathrm{c}}
\def\rd{\mathrm{d}}
\def\ve{\mathrm{ve}}
\def\ef{\mathrm{ef}}
\def\VE{\mathrm{VE}}
\def\EF{\mathrm{EF}}
\def\a{\alpha}
\def\b{\beta}
\def\ga{\gamma}
\def\d{\delta}  
\def\n{\mbox{\boldmath $n$}}
\title{\bf Triakis Solids and Harmonic Functions} 
\author{Miki Hasui and Katsunori Iwasaki\thanks{Department of Mathematics, 
Hokkaido University, Kita 10, Nishi 8, Kita-ku, Sapporo 060-0810 Japan. 
E-mail: {\tt s113027@math.sci.hokudai.ac.jp} and 
{\tt iwasaki@math.sci.hokudai.ac.jp} (corresponding author)}}
\date{October 17, 2012} 
\begin{document}
\maketitle
\begin{abstract} 
We describe the harmonic functions for certain isohedral triakis solids.   
They are the first examples for which polyhedral harmonics is strictly 
larger than group harmonics. \\[3mm]
\end{abstract} 
\section{Introduction} \label{sec:intro} 
A harmonic function on $\R^n$ is a distribution solution to the Laplace 
equation 
\begin{equation} \label{eqn:laplacian}
\varDelta f = (\partial_1^2 + \cdots + \partial_n^2)f = 0, \qquad 
\partial_i = \partial/\partial x_i, 
\end{equation}
which is necessarily a smooth function.  
A classical theorem of Gauss and Koebe states that a continuous 
function on $\R^n$ is harmonic if and only if it satisfies the mean 
value property with respect to the sphere $S^{n-1}$.  
On the other hand the Laplacian $\varDelta$ admits an invariant-theoretic 
interpretation relative to the orthogonal group $O(n)$. 
Namely the $\R$-algebra of $O(n)$-invariant polynomials in 
$x = (x_1, \dots, x_n)$ is generated by the squared distance function 
$\varphi(x) = x_1^2 + \cdots + x_n^2$ and the Laplacian $\varDelta = 
\varphi(\partial)$ is the result of substituting $\partial = 
(\partial_1, \dots, \partial_n)$ into $\varphi(x)$.  
These two characterizations of harmonic functions allow us to generalize 
the notion in two directions, that is, to polyhedral harmonics and to 
group harmonics.  
\par
Let $P$ be an $n$-dimensional polytope in $\mathbb{R}^n$ and 
$P(k)$ the $k$-dimensional skeleton of $P$, that is, the union 
of all $k$-dimensional faces of $P$. 
A continuous function $f(x)$ on $\R^n$ is said to be $P(k)$-{\sl harmonic} 
if it satisfies the mean value property with respect to $P(k)$, that is, if
\[
f(x) = \dfrac{1}{|P(k)|} \int_{P(k)} f(x+ r y) \, d\mu_k(y) \qquad 
(\, {}^{\forall} x \in \R^n, \, {}^{\forall} r > 0),  
\]
where $\mu_k$ is the $k$-dimensional Euclidean measure on $P(k)$ with 
$|P(k)| := \mu_k(P(k))$ being its total volume. 
Let $\H_{P(k)}$ denote the set of all $P(k)$-harmonic functions on $\R^n$. 
A general result of Iwasaki \cite[Theorem 1.1]{Iwasaki1} states that 
$\H_{P(k)}$ is a finite-dimensional linear space of polynomials. 
It is naturally an $\R[\partial]$-module, that is, stable under partial 
differentiations. 
If $G(k)$ is the symmetry group of $P(k)$, then it is also  
an $\R[G(k)]$-module, that is, stable under the natural action of $G(k)$.  
Moreover, if $P(k)$ enjoys ample symmetry, that is, if $G(k)$ acts on 
$\R^n$ irreducibly then $\H_{P(k)}$ is a finite-dimensional linear space 
of harmonic polynomials \cite[Theorem 2.4]{Iwasaki1}.   
\par
Let $G$ be a subgroup of $O(n)$ and $\R[x]^G$ the $\R$-algebra of 
$G$-invariant polynomials of $x = (x_1,\dots,x_n)$. 
A $G$-{\sl harmonic} function on $\R^n$ is a distribution solution 
to the system of PDEs:  
\begin{equation}  \label{eqn:g-harmonic}
\varphi(\partial) f = 0 \qquad (\, {}^{\forall} \varphi(x) \in 
\R[x]^G_+), 
\end{equation} 
where $\R[x]^G_+$ is the maximal ideal of $\R[x]^G$ consisting of 
those $\varphi(x)$'s without constant term: $\varphi(0) = 0$ (see 
Helgason \cite[Chap. I\!I\!I]{Helgason}). 
Let $\H_G$ denote the set of all $G$-harmonic functions on $\R^n$. 
To define system \eqref{eqn:g-harmonic}, the polynomial $\varphi(x)$ 
may not range over all $\R[x]^G_+$, but only over a set of generators 
of $\R[x]^G_+$.    
For example, if $G$ is the entire $O(n)$ then $\R[x]^{O(n)}_+$ is 
generated by the squared distance function only, so that system 
\eqref{eqn:g-harmonic} reduces to the single Laplace equation 
\eqref{eqn:laplacian}. 
In this article, $G$ will be a finite group in general and a finite 
reflection group in particular. 
Steinberg \cite{Steinberg} shows that if $G$ is a finite group 
then $\H_G$ is a finite-dimensional linear space of polynomials 
with $\dim \H_G \ge |G|$, the order of $G$. 
When $G$ is a finite reflection group, he goes on to determine 
$\H_G$ explicitly. 
Namely, as an $\R[\partial]$-module $\H_G$ is generated by the 
fundamental alternating polynomial $\varDelta_G(x)$ of $G$, 
and as an $\R[G]$-module it is the regular representation of 
$G$, in particular $\dim \H_G = |G|$.  
\par
When $G = G(k)$ is the symmetry group of $P(k)$, it makes sense 
to compare $\H_{P(k)}$ with $\H_{G(k)}$. 
According to a result in \cite[formula (2.14)]{Iwasaki1}  
there is always the inclusion 
\begin{equation} \label{eqn:inclusion}
\H_{G(k)} \subset \H_{P(k)} \qquad (k = 0, \dots, n). 
\end{equation} 
It is known that $\H_{P(k)}$ coincides with $\H_{G(k)}$ when $P$ is any 
regular convex polytope with center at the origin, in which case 
$G(k)$ is an irreducible finite reflection group so that $\H_{P(k)} = 
\H_{G(k)}$ is determined by Steinberg's theorem; see 
Iwasaki \cite[Theorem 4.4]{Iwasaki4} and the references therein. 
One can show that the coincidence also occurs, for example, for the 
truncated icosahedron $[5,6,6]$ of an Archimedean solid, or the soccer ball.     
So far, however, no polytope with ample symmetry has been known for 
which inclusion \eqref{eqn:inclusion} is strict (for some $k > 0$), 
that is, polyhedral harmonics is strictly larger than group harmonics. 
The aim of this article is to present the first examples of such polytopes. 
They arise from a one-parameter family of isohedral triakis tetrahedra 
as well as from a one-parameter family of isohedral triakis octahedra 
in three dimensions. 
The former family contains a desired example for $k = 1$, but none for 
$k = 0, 2, 3$ (see Theorem \ref{thm:main-tetra}), while the latter 
family contains such an example for every $k = 0,1,2,3$ (see 
Theorem \ref{thm:main-octa}). 
In this respect the latter family is more interesting than the former, 
but in any case we begin with the simpler case of the former family and 
then proceed to the latter.       
\section{A Family of Isohedral Triakis Tetrahedra} \label{sec:family-tetra}
\begin{figure}[t]
\begin{center}
\unitlength 0.1in
\begin{picture}( 29.0000, 26.8000)(  3.9000,-28.5000)
%
\special{pn 20}%
\special{pa 1790 400}%
\special{pa 2380 2810}%
\special{fp}%
%
\special{pn 20}%
\special{pa 3180 2190}%
\special{pa 2380 2800}%
\special{fp}%
%
\special{pn 20}%
\special{pa 2370 2810}%
\special{pa 590 2190}%
\special{fp}%
%
\special{pn 20}%
\special{pa 1790 420}%
\special{pa 610 2180}%
\special{fp}%
%
\special{pn 13}%
\special{pa 600 2190}%
\special{pa 2120 2200}%
\special{da 0.070}%
%
\special{pn 20}%
\special{pa 1800 420}%
\special{pa 3180 1060}%
\special{fp}%
%
\special{pn 20}%
\special{pa 3180 1060}%
\special{pa 3180 2200}%
\special{fp}%
%
\special{pn 20}%
\special{pa 3180 1060}%
\special{pa 2380 2810}%
\special{fp}%
%
\special{pn 13}%
\special{pa 2920 1860}%
\special{pa 3190 2200}%
\special{da 0.070}%
%
\special{pn 13}%
\special{pa 2320 2200}%
\special{pa 2590 2200}%
\special{da 0.070}%
%
\special{pn 13}%
\special{pa 2770 2200}%
\special{pa 3180 2200}%
\special{da 0.070}%
\put(17.5000,-3.4000){\makebox(0,0)[lb]{$\rA$}}%
\put(23.6000,-30.2000){\makebox(0,0)[lb]{$\rB$}}%
\put(32.9000,-22.9000){\makebox(0,0)[lb]{$\rC$}}%
\put(3.9000,-22.7000){\makebox(0,0)[lb]{$\rD$}}%
%
\special{pn 13}%
\special{pa 3180 1070}%
\special{pa 2180 1620}%
\special{fp}%
%
\special{pn 13}%
\special{pa 2030 1700}%
\special{pa 1790 1820}%
\special{fp}%
%
\special{pn 13}%
\special{pa 1800 440}%
\special{pa 2540 1360}%
\special{da 0.070}%
%
\special{pn 13}%
\special{pa 2620 1470}%
\special{pa 2810 1700}%
\special{da 0.070}%
%
\special{pn 20}%
\special{sh 1.000}%
\special{ar 2310 1550 42 42  0.0000000 6.2831853}%
%
\special{pn 20}%
\special{sh 1.000}%
\special{ar 1800 1830 42 42  0.0000000 6.2831853}%
\put(16.0000,-19.9000){\makebox(0,0)[lb]{$\rO$}}%
\put(32.7000,-10.7000){\makebox(0,0)[lb]{$\rd$}}%
\put(22.7000,-14.4000){\makebox(0,0)[lb]{$\rd'$}}%
\end{picture}%
\end{center}
\caption{Adjoining a pyramid to a face of a regular tetrahedron; 
$\overline{\rO\rd} : \overline{\rO\rd'} = r : 1$.}
\label{fig:triakis-tetra}
\end{figure}
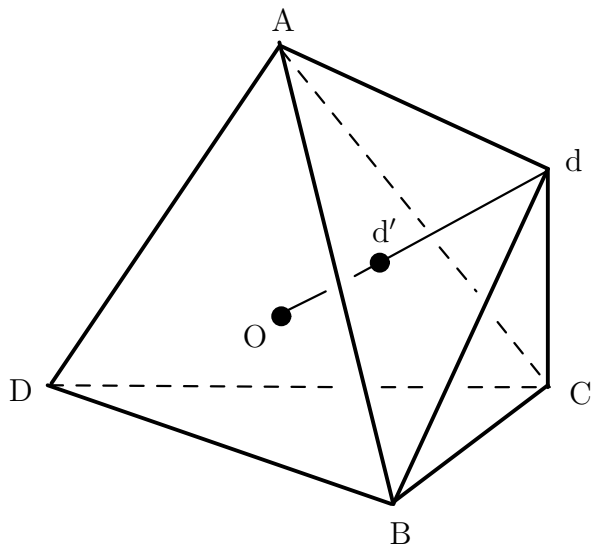
An {\sl isohedral triakis tetrahedron} is obtained from a regular 
tetrahedron by adjoining to each face of it a pyramid (based on that 
face) of appropriate height, or excavating such a pyramid, where 
a polyhedron is said to be {\sl isohedral} if its symmetry group is 
transitive on its faces, that is, if all faces are equivalent under 
symmetries of the polyhedron; see Gr\"{u}nbaum and Shephard \cite{GS}.     
Our triakis tetrahedron $P$ depends on a positive parameter $r$.  
To describe it more neatly, let $T$ be the regular tetrahedron  
with which we get started; it is centered at the origin $\rO$ and 
with four vertices $\rA$, $\rB$, $\rC$, $\rD$. 
If $\rd'$ denotes the center of the face $\rA\rB\rC$, then the 
pyramid based on this face has the top vertex $\rd$ that lies on 
the ray emanating from $\rO$ and passing through $\rd'$ in such 
a manner that the distance ratio 
$\overline{\rO\rd} : \overline{\rO\rd'}$ is $r : 1$ 
(see Figure \ref{fig:triakis-tetra}).    
The polyhedron $P$ has four six-valent vertices  
$\rA$, $\rB$, $\rC$, $\rD$ and four three-valent vertices  
$\ra$, $\rb$, $\rc$, $\rd$, where $\ra$, $\rb$, $\rc$ are defined 
in a similar manner as the vertex $\rd$ is, and $P$ has twelve 
faces and eighteen edges.  
\par
\begin{figure}[p]
\begin{minipage}{0.49\hsize}
\begin{center}
\includegraphics*[width=60mm,height=60mm]{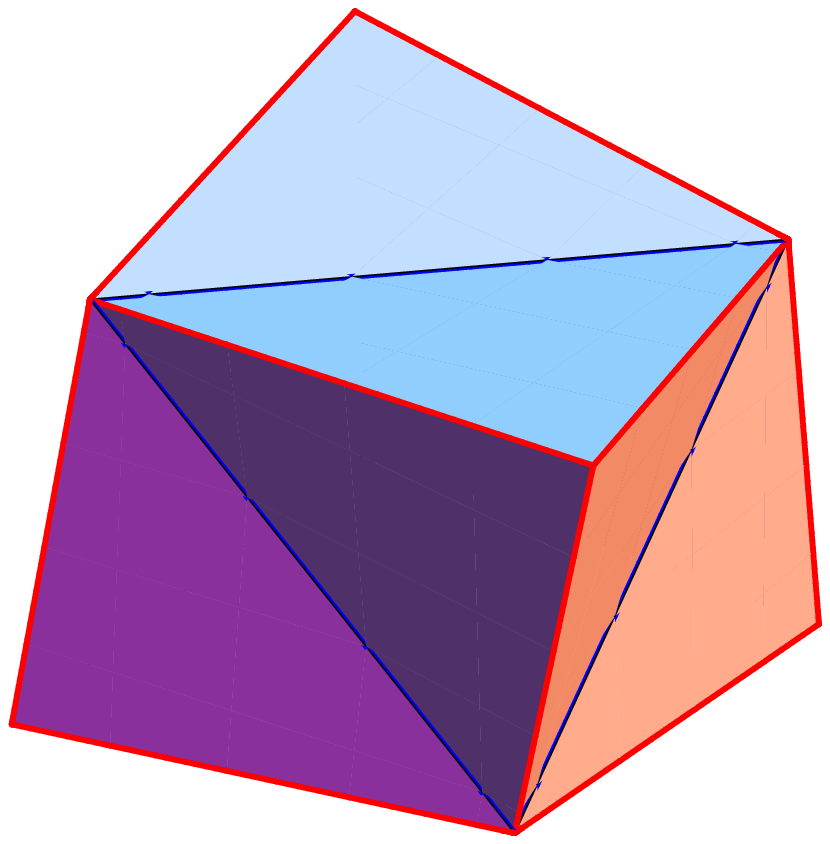} 
\end{center}
\vspace{-7mm} 
\centerline{$r > 3$}
\end{minipage}
\hspace{-5mm}
\begin{minipage}{0.49\hsize}
\begin{center}
\includegraphics*[width=60mm,height=60mm]{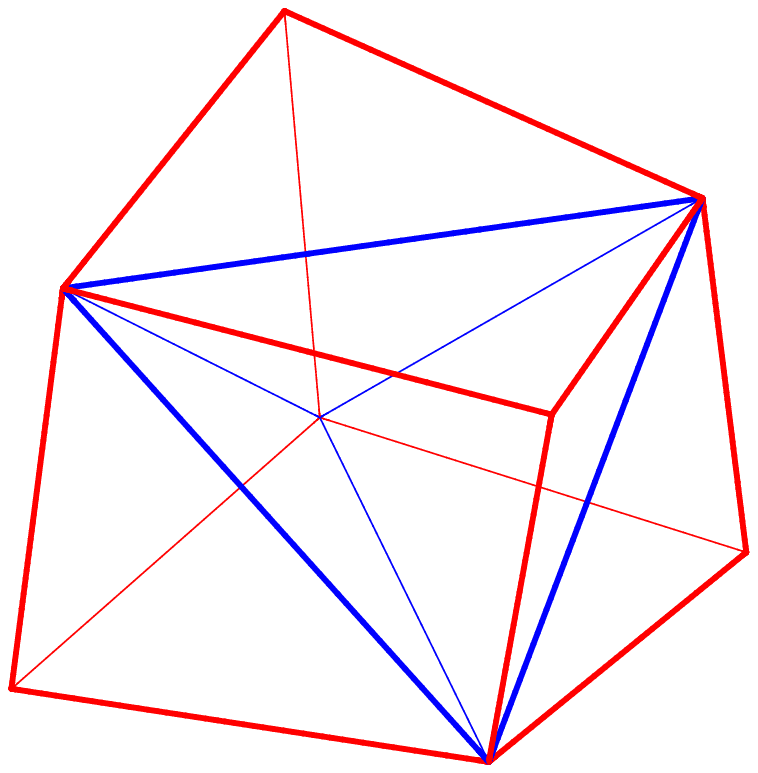} 
\end{center}
\vspace{-7mm}
\centerline{$r > 3$}
\end{minipage}
\par\vspace{8mm}\noindent
\begin{minipage}{0.49\hsize}
\begin{center}
\includegraphics*[width=55mm,height=55mm]{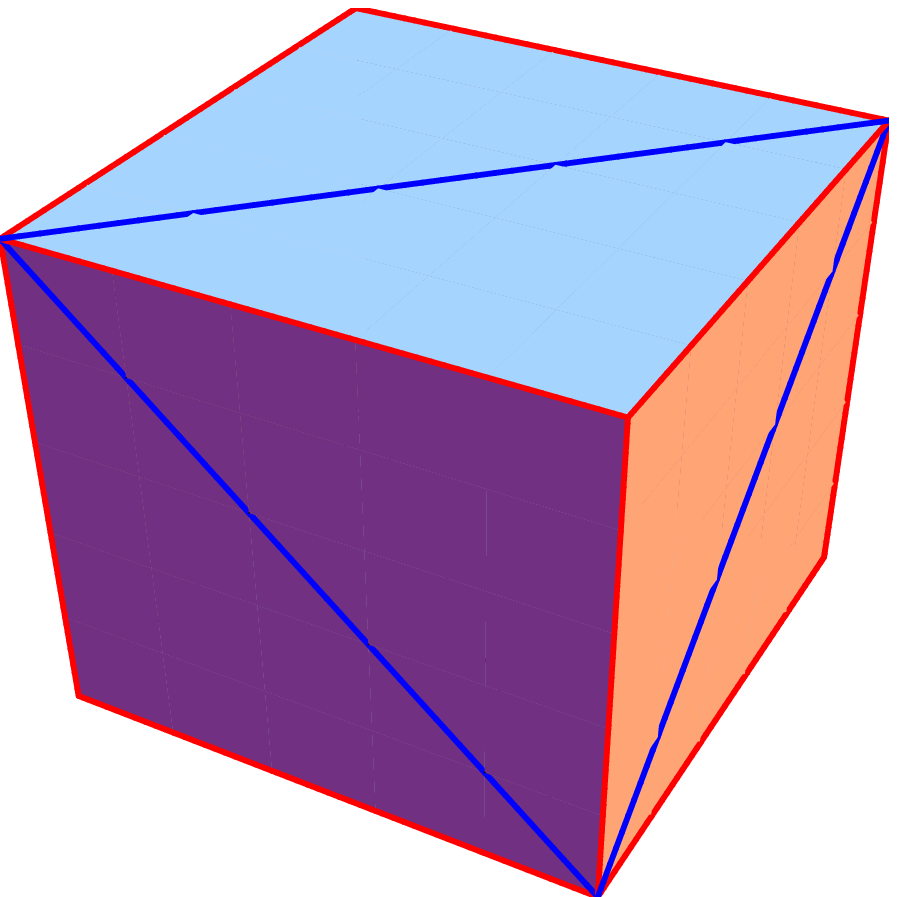} 
\end{center} 
\centerline{$r = 3$}
\end{minipage}
\hspace{-5mm}
\begin{minipage}{0.49\hsize}
\begin{center}
\includegraphics*[width=55mm,height=55mm]{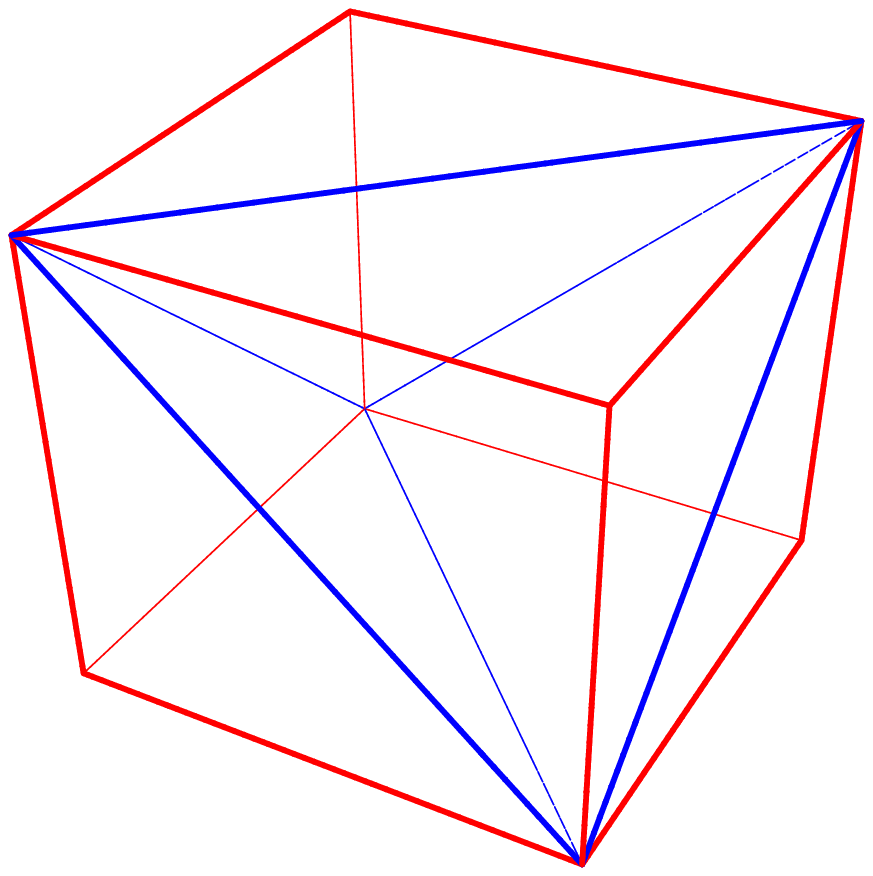} 
\end{center} 
\centerline{$r = 3$}
\end{minipage}
\par\vspace{8mm}\noindent
\begin{minipage}{0.49\hsize}
\begin{center}
\includegraphics*[width=60mm,height=60mm]{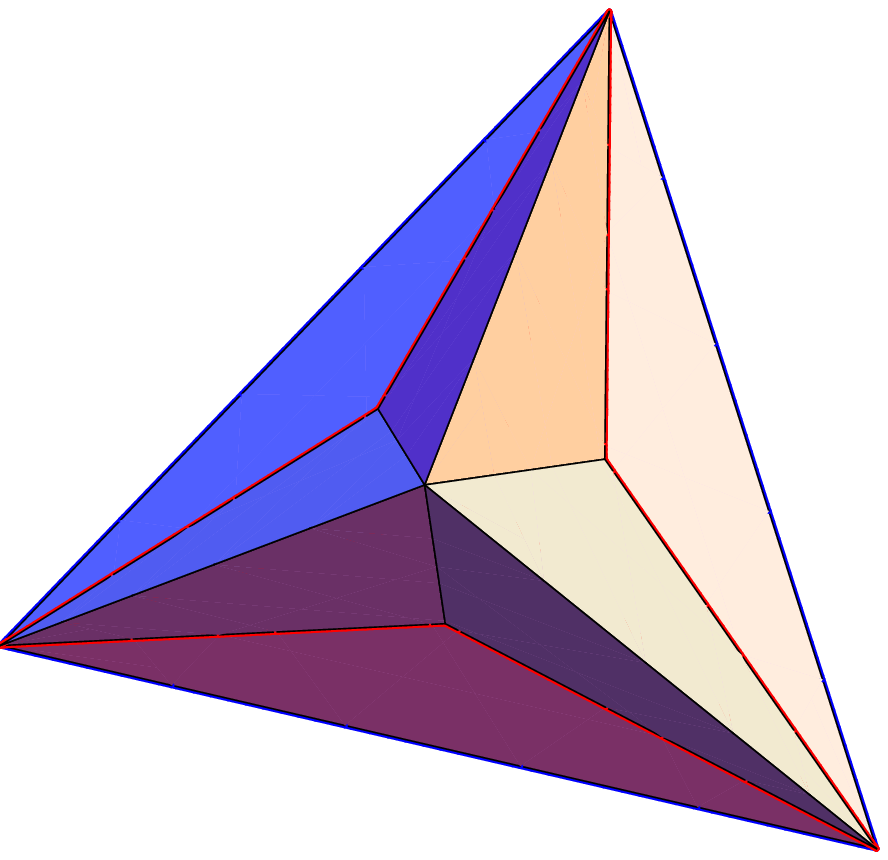} 
\end{center}
\vspace{-10mm} 
\centerline{$0 < r < 1$}
\end{minipage}
\hspace{-5mm}
\begin{minipage}{0.49\hsize}
\begin{center}
\includegraphics*[width=60mm,height=60mm]{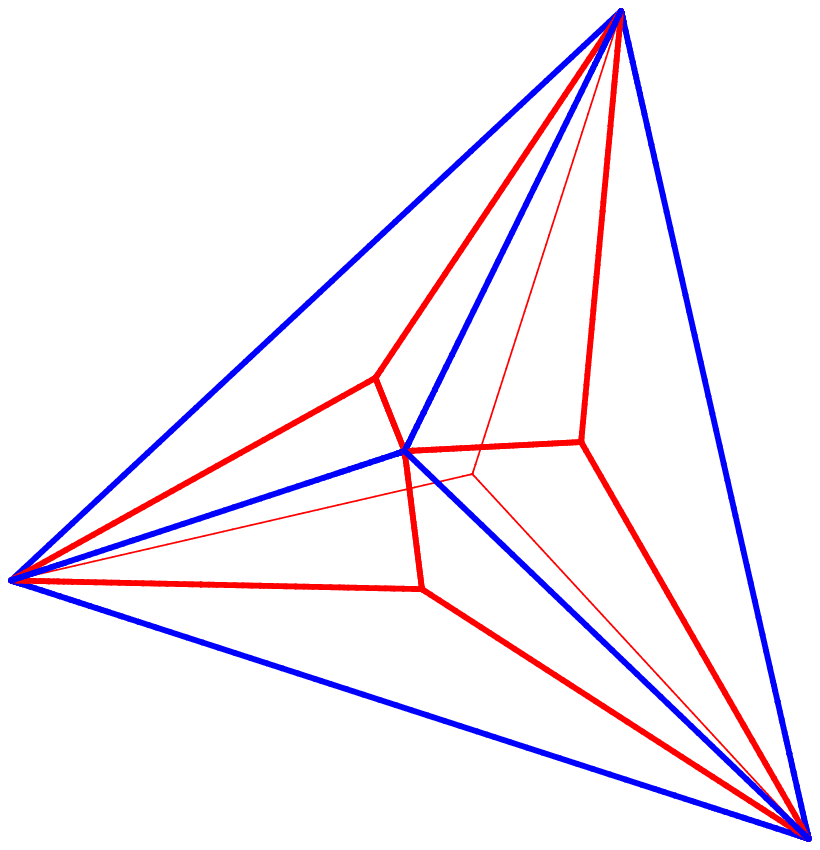} 
\end{center}
\vspace{-10mm} 
\centerline{$0 < r < 1$}
\end{minipage}
\par\vspace{8mm}\noindent
\caption{A one-parameter family of isohedral triakis tetrahedra. 
Each polyhedron is shown in a ``cardboard" view on the left, as 
well as in a skeletal view on the right. Polyhedron with $1 < r < 3$ 
is not indicated, but imaginable from the one with $r > 3$, although 
the former is convex while the latter is not. See also a series of 
pictures in Gr\"{u}nbaum \cite[Fig. 14]{Grunbaum}. }
\label{fig:family-tetra}
\end{figure}

Let us look more closely at the polyhedron $P$ for various values 
of $r$ (see Figure \ref{fig:family-tetra}). 
The values $r = 1$ and $r = 3$ are special in that as a point set, $P$ 
degenerates to the tetrahedron $T$ at $r = 1$ and it becomes a 
cube $C$ at $r = 3$.  
For $r > 1$, $P$ is obtained from $T$ by adjoining a pyramid to each 
face of $T$, while for $0 < r < 1$, it is obtained by excavating 
such a pyramid.   
Moreover $P$ is convex if and only if $1 \le r \le 3$, in which interval 
the value $r = 9/5$ is distinguished in that $P$ becomes a Catalan 
solid \cite{Catalan,KKK}, or an Archimedean dual solid, whose dual is 
the truncated tetrahedron $[3,6,6]$. 
When $r = 1$, although the four points $\rA$, $\rB$, $\rC$, $\rd$ are 
coplanar, we think of $\rA\rB\rd$, $\rB\rC\rd$, $\rA\rC\rd$ as   
distinct faces of $P$, and thus $\rA\rd$, $\rB\rd$, $\rC\rd$ as  
edges, and $\rd$ as a vertex of $P$. 
Similarly, when $r = 3$, although the four points 
$\rA$, $\rB$, $\rc$, $\rd$ are coplanar, we think of $\rA\rB\rc$ and  
$\rA\rB\rd$ as distinct faces of $P$, and thus $\rA\rB$ as 
an edge of $P$. 
With this convention, as a combinatorial polyhedron, $P$ has the 
constant skeletal structure for all $r > 0$.     
\par
Let $G(k)$ be the symmetry group of $P(k)$. 
When $k = 1$, for every $r > 0$ the group $G(k)$ is always the same as 
the symmetry group of the tetrahedron $T$, which we denote by $W(\rA_3)$ 
since it is a Weyl group of type $\rA_3$. 
On the other hand, when $k = 0, 2, 3$, the group $G(k)$ stays the same as 
$W(\rA_3)$ as long as $r \neq 3$, but at $r = 3$ it jumps up to be the 
symmetry group of the cube $C$, which is denoted by $W(\rB_3)$ as being a 
Weyl group of type $\rB_3$.  
For simplicity of notation, $\H_{W(\rA_3)}$ and $\H_{W(\rB_3)}$ are  
abbreviated to $\H_{\rA_3}$ and $\H_{\rB_3}$ respectively. 
Note that
\[
\H_{\rA_3} \varsubsetneq \H_{\rB_3}, \qquad 
\dim \H_{\rA_3} = 24, \qquad \dim \H_{\rB_3} = 48. 
\]
\par
It is our problem to consider how $\H_{P(k)}$ behaves as $r$  
varies and when $\H_{P(k)}$ is strictly larger than $\H_{G(k)}$.  
Any value of $r$ with this phenomenon is referred 
to as a {\sl critical value}. 
It turns out that the vertex $(k = 0)$ and face $(k =2)$ 
problems have no critical values, but the edge $(k = 1)$ 
problem certainly has a critical value $r_1 = 3.62398\cdots$. 
It is an algebraic integer of degree six that is the unique 
positive root of a sextic equation  
\begin{equation} \label{eqn:sextic}
\chi_1(r) := r^6 + 2 r^5 + r^4 - 36 r^3  - 45 r^2 - 270 r - 405 = 0.  
\end{equation}
How this algebraic equations arises will be explained in 
\S\ref{subsec:edge-tetra}.   
The volume problem $(k = 3)$ need not be dealt with independently,  
since the isohedrality of $P$ implies $\H_{P(2)} = \H_{P(3)}$ by a result 
of \cite[Theorem 2.2]{Iwasaki1} so that the face and volume problems 
have the same solution.   
\begin{theorem} \label{thm:main-tetra} 
For the family of isohedral triakis tetrahedra, $\H_{P(k)}$ is strictly 
larger than $\H_{G(k)}$ if and only if $k = 1$ and $r = r_1 = 3.62398\cdots$.  
The function space $\H_{P(k)}$ is given by  
\begin{equation} \label{eqn:main-tetra}
\H_{P(k)} = \begin{cases}
\H_{\rB_3} & (\mbox{if either $k = 1$, $r = r_1;$ or $k = 0, 2, 3$, $r = 3$}), \\ 
\H_{\rA_3} & (\mbox{otherwise}). 
\end{cases} 
\end{equation}  
\end{theorem}
\par
This theorem asserts that when $k = 1$ and $r$ is at the critical 
value $r_1$, the figure $P(1)$ has only the same symmetries as a tetrahedron, 
but $\H_{P(1)}$ jumps up to become the space of cubic harmonics, which is 
strictly larger than that of tetrahedral harmonics. 
For each $k = 0, 2, 3$, a jumping phenomenon also occurs at $r = 3$ 
with the space $\H_{P(k)}$ jumping from $\H_{\rA_3}$ to $\H_{\rB_3}$, 
but at the same time the group $G(k)$ also jumps from $W(\rA_3)$ to 
$W(\rB_3)$. 
Altogether, the equality $\H_{P(k)} = \H_{G(k)}$ continues 
to hold and no critical phenomenon occurs at $r = 3$. 
After a review in \S\ref{sec:inv-pde} on PDEs that characterize polyhedral 
harmonics, Theorem \ref{thm:main-tetra} will be established 
in \S\ref{sec:skeleton-tetra}.          
\section{Invariant Differential Equations} \label{sec:inv-pde}
Iwasaki \cite[Theorem 2.1]{Iwasaki1} derives a system of partial 
differential equations  
\begin{equation} \label{eqn:pde}
\tau_m^{(k)}(\partial) f = 0 \qquad (m = 1,2,3,\dots)  
\end{equation}
that characterizes $\H_{P(k)}$ as its distribution solution space.  
Here is a brief review of it when $P$ is a three-dimensional 
polyhedron and $k = 0, 1, 2$.      
For $j = 0,1,2$, let $\{P_{i_j}\}_{i_j \in I_j}$ be the set of all 
$j$-dimensional faces of $P$, where $I_j$ is an index set.   
Let $H_{i_j}$ be the $j$-dimensional affine subspace of $\R^3$ 
containing $P_{i_j}$. 
Let $p_{i_j}$ be the foot of orthogonal projection from the 
origin $\rO$ down to $H_{i_j}$. 
We mean by $i_j \prec i_{j+1}$ that $P_{i_j}$ is a face of 
$P_{i_{j+1}}$. 
If $i_j \prec i_{i_{j+1}}$ then the vector $p_{i_j} - p_{i_{j+1}}$ 
is parallel to the outer unit normal vector $\n_{i_j,i_{j+1}}$ of 
$\partial P_{i_{j+1}}$ in $H_{i_{j+1}}$ at the face 
$P_{i_j}$, so that a number $[i_j:i_{j+1}]$, called the 
{\sl incidence number}, is defined by the relation 
$p_{i_j} - p_{i_{j+1}} = [i_j:i_{j+1}] \n_{i_j,i_{j+1}}$. 
Put $I(1) = \{ i = (i_0, i_1) \,:\, i_0 \prec i_1\}$,  
$I(2) =  \{ i = (i_0, i_1, i_2) \,:\, i_0 \prec i_1 \prec i_2\}$, 
and define 
\[
[i] = \begin{cases}
[i_0:i_1]  & (i = (i_0,i_1) \in I(1)), \\ 
[i_0:i_1][i_1:i_2] & (i = (i_0,i_1,i_2) \in I(2)). 
\end{cases}
\]
Then $\tau_m^{(k)}(x)$, $k = 0,1,2$, are homogeneous 
polynomials of degree $m$ defined by 
\begin{align}
\tau_m^{(0)}(x) &:= \sum_{i_0 \in I_0} \langle p_{i_0}, x \rangle^m, 
\label{eqn:tau0} \\
\tau_m^{(1)}(x) &:= \sum_{i \in I(1)} [i] \, 
h_m(\langle p_{i_0}, x \rangle, \, \langle p_{i_1}, x\rangle), 
\label{eqn:tau1} \\
\tau_m^{(2)}(x) &:= \sum_{i \in I(2)} [i] \, 
h_m(\langle p_{i_0}, x \rangle, \, \langle p_{i_1}, x \rangle, \, 
\langle p_{i_2}, x \rangle), \label{eqn:tau2} 
\end{align}  
where $\langle \,\cdot\, , \, \cdot \, \rangle$ is the standard 
inner product on $\R^3$ and $h_m$ stands for the complete 
symmetric polynomial of degree $m$ in two or three 
variables (see e.g. Macdonald \cite{Macdonald}).   
\par
\begin{figure}[t]
\begin{center}
\unitlength 0.1in
\begin{picture}( 45.5000, 34.8000)( 14.1000,-40.8000)
%
\special{pn 20}%
\special{pa 3220 1330}%
\special{pa 3810 3740}%
\special{fp}%
%
\special{pn 20}%
\special{pa 3230 1350}%
\special{pa 4610 1990}%
\special{fp}%
%
\special{pn 20}%
\special{pa 4610 1990}%
\special{pa 3810 3740}%
\special{fp}%
\put(32.1000,-12.8000){\makebox(0,0)[lb]{$\rA$}}%
\put(38.0000,-39.1000){\makebox(0,0)[lb]{$\rB$}}%
%
\special{pn 13}%
\special{pa 3460 2630}%
\special{pa 3220 2750}%
\special{fp}%
%
\special{pn 20}%
\special{sh 1.000}%
\special{ar 3900 2390 42 42  0.0000000 6.2831853}%
%
\special{pn 20}%
\special{sh 1.000}%
\special{ar 3230 2760 42 42  0.0000000 6.2831853}%
\put(30.3000,-29.2000){\makebox(0,0)[lb]{$\rO$}}%
\put(46.2000,-19.4000){\makebox(0,0)[lb]{$\rd$}}%
%
\special{pn 13}%
\special{pa 3410 610}%
\special{pa 5960 2120}%
\special{fp}%
%
\special{pn 13}%
\special{pa 5960 2110}%
\special{pa 3950 4080}%
\special{fp}%
%
\special{pn 13}%
\special{pa 3420 600}%
\special{pa 1410 2570}%
\special{fp}%
%
\special{pn 13}%
\special{pa 1410 2550}%
\special{pa 3960 4070}%
\special{fp}%
%
\special{pn 13}%
\special{pa 3620 2540}%
\special{pa 3910 2390}%
\special{fp}%
%
\special{pn 13}%
\special{pa 2650 1080}%
\special{pa 5710 2500}%
\special{fp}%
\put(20.6000,-24.0000){\makebox(0,0)[lb]{$H_{\rA\rB\rd}$}}%
\put(23.5000,-10.6000){\makebox(0,0)[lb]{$\ell_{\rA\rd}$}}%
%
\special{pn 20}%
\special{sh 1.000}%
\special{ar 4050 1730 42 42  0.0000000 6.2831853}%
%
\special{pn 13}%
\special{pa 3230 2760}%
\special{pa 3440 2480}%
\special{fp}%
%
\special{pn 13}%
\special{pa 3530 2360}%
\special{pa 4050 1720}%
\special{fp}%
\put(39.0000,-26.1000){\makebox(0,0)[lb]{$\rA\rB\rd$}}%
%
\special{pn 13}%
\special{pa 4000 1990}%
\special{pa 4120 1470}%
\special{fp}%
\special{sh 1}%
\special{pa 4120 1470}%
\special{pa 4086 1530}%
\special{pa 4108 1522}%
\special{pa 4124 1540}%
\special{pa 4120 1470}%
\special{fp}%
\put(41.3000,-17.2000){\makebox(0,0)[lb]{$\rA\rd$}}%
%
\special{pn 13}%
\special{pa 4040 1830}%
\special{pa 4120 1870}%
\special{pa 4140 1780}%
\special{pa 4140 1780}%
\special{pa 4140 1780}%
\special{fp}%
\put(40.7000,-14.5000){\makebox(0,0)[lb]{$\n_{\rA\rd, \rA\rB\rd}$}}%
%
\special{pn 13}%
\special{pa 3130 1310}%
\special{pa 3020 1250}%
\special{fp}%
\special{sh 1}%
\special{pa 3020 1250}%
\special{pa 3070 1300}%
\special{pa 3068 1276}%
\special{pa 3088 1264}%
\special{pa 3020 1250}%
\special{fp}%
\put(27.7000,-14.6000){\makebox(0,0)[lb]{$\n_{\rA, \rA\rd}$}}%
%
\special{pn 13}%
\special{pa 4690 2030}%
\special{pa 4850 2100}%
\special{fp}%
\special{sh 1}%
\special{pa 4850 2100}%
\special{pa 4798 2056}%
\special{pa 4802 2080}%
\special{pa 4782 2092}%
\special{pa 4850 2100}%
\special{fp}%
\put(48.1000,-20.8000){\makebox(0,0)[lb]{$\n_{\rd, \rA\rd}$}}%
\end{picture}%
\end{center}
\caption{Combinatorial data for the system \eqref{eqn:pde} of PDEs.}
\label{fig:foot}
\end{figure}
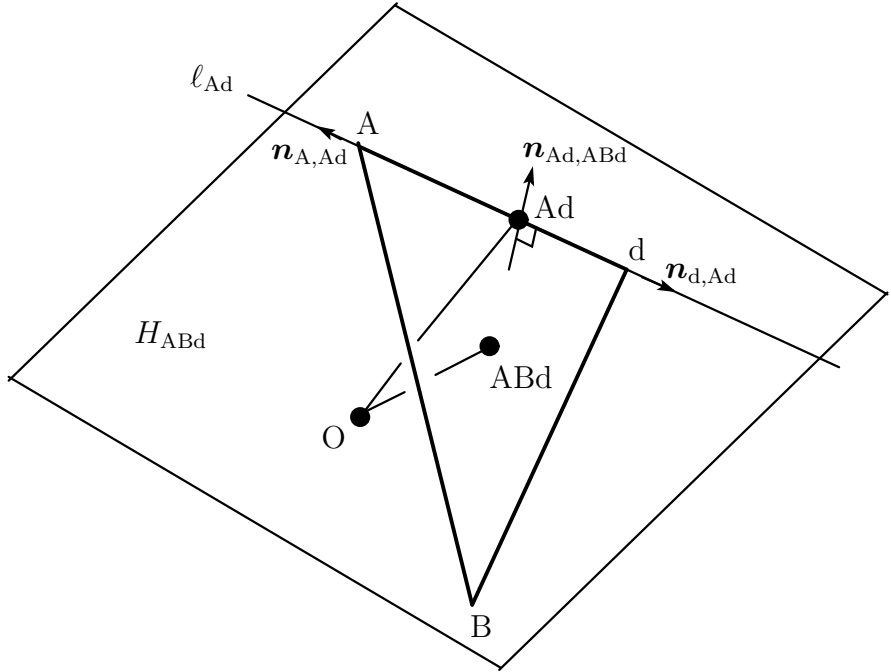

The general construction mentioned above will be applied to the 
particular case of our polyhedron $P$ upon adjusting the notation 
to the current situation. 
In order to represent the index sets $I_0$, $I_1$, $I_2$, it would 
be best to let the vertices, edges and faces to speak of themselves:    
\begin{align*}
I_0 &= \{\rA, \rB, \rC, \rD \} \cup \{\ra, \rb, \rc, \rd\},  \\
I_1 &= \{\rA\rB, \rA\rC, \rA\rD, \rB\rC, \rB\rD, \rC\rD \} \cup 
\{\rA\rb, \rA\rc, \rA\rd, \rB\ra, \rB\rc, \rB\rd,  
\rC\ra, \rC\rb, \rC\rd, \rD\ra, \rD\rb, \rD\rc\}, \\
I_2 &= \{\rA\rB\rc, \rA\rB\rd, \rA\rC\rb, \rA\rC\rd, \rA\rD\rb, \rA\rD\rc,   
\rB\rC\ra, \rB\rC\rd, \rB\rD\ra, \rB\rD\rc, \rC\rD\ra, \rC\rD\rb\},   
\end{align*}
where $\{\cdots\}$ stands for an orbit under symmetries. 
For an index $\rA\rd \in I_1$, the same symbol $\rA\rd$ denotes 
the foot of orthogonal projection from the origin $\rO$ to the affine 
line $\ell_{\rA\rd}$ passing through $\rA$ and $\rd$; this rule also 
applies to another index $\rA\rB \in I_1$ as well as to all the other 
indices of $I_1$. 
In a similar manner, for an index $\rA\rB\rd \in I_2$, the same 
symbol $\rA\rB\rd$ denotes the foot of orthogonal projection from 
$\rO$ to the affine plane $H_{\rA\rB\rd}$ passing through 
$\rA$, $\rB$, $\rd$, with this rule applying to all the other 
indices of $I_2$; see Figure \ref{fig:foot}.    
\par
Up to symmetries there are three types of vertex-to-edge 
incidence relations: 
\begin{equation} \label{eqn:VE-tetra}
(\VE 1) \quad \rA \prec \rA\rB, \qquad 
(\VE 2) \quad \rA \prec \rA\rd, \qquad 
(\VE 3) \quad \rd \prec \rA\rd.  
\end{equation}
On the other hand, up to symmetries there are two types of 
edge-to-face incidence relations:  
\begin{equation} \label{eqn:EF-tetra}
(\EF 1) \quad \rA\rB \prec \rA\rB\rd, \qquad 
(\EF 2) \quad \rA\rd \prec \rA\rB\rd.  
\end{equation}
The incidence number of an incidence relation depends only on 
its type. 
Let $\ve_{\nu}$ denote the incidence number of type ($\VE \nu$), 
$\nu = 1,2,3$, and $\ef_{\nu}$ denote that of type ($\EF \nu$), 
$\nu = 1,2$.    
For example we have $[\rA:\rA\rd] = \ve_2$ and    
$[\rA\rd : \rA\rB\rd] = \ef_2$, so that  
$\rA - \rA\rd = \ve_2 \, \n_{\rA, \rA\rd}$ and 
$\rA\rd - \rA\rB\rd = \ef_2 \, \n_{\rA\rd, \rA\rB\rd}$  
for $\rA \prec \rA\rd \prec \rA\rB\rd$. 
These numbers are positive in the situation of Figure \ref{fig:foot}, 
but they or some other incidence numbers may be negative for some 
values of $r$.   
\par
For explicit calculations, it is convenient to work with 
Cartesian coordinates by putting 
\begin{equation} \label{eqn:coord-tetra}
\rA = (1,-1,-1), \quad \rB = (-1,1,-1), \quad \rC = (-1,-1,1), 
\quad \rD = (1,1,1). 
\end{equation}
The symmetry group $W(\rA_3)$ of the tetrahedron $T$ then admits an 
invariant basis 
\[
e_2(x) := x_1^2 + x_2^2 + x_3^2, \quad 
e_3(x) := x_1 x_2 x_3, \quad 
e_4(x) := x_2^2 x_3^2 + x_3^2 x_1^2 + x_1^2 x_2^2,    
\]
so that $\H_{\rA_3}$ can be characterized as the solution 
space to the system of PDEs: 
\begin{equation} \label{eqn:pde-a3}
e_2(\partial) f = e_3(\partial) f = e_4(\partial) f = 0.
\end{equation}
As an $\R[\partial]$-module, $\H_{\rA_3}$ is generated by the 
fundamental alternating polynomial 
\[
\varDelta_{\rA_3}(x) = (x_1^2-x_2^2)(x_2^2-x_3^2)(x_3^2-x_1^2).  
\]
\par
If the vertices of $T$ are taken as in formula \eqref{eqn:coord-tetra}, 
then $C$ is the cube with vertices at $(\pm1, \pm1, \pm 1)$. 
The symmetry group $W(\rB_3)$ of $C$ then admits an invariant basis 
$e_2(x)$, $e_4(x)$, $e_6(x) := e_3^2(x)$, so that $\H_{\rB_3}$ can 
be characterized as the solution space to the system of PDEs: 
\begin{equation} \label{eqn:pde-b3}
e_2(\partial) f = e_4(\partial) f = e_6(\partial) f = 0.
\end{equation}
As an $\R[\partial]$-module, $\H_{\rB_3}$ is generated by the 
fundamental alternating polynomial 
\begin{equation} \label{eqn:fap-b3}
\varDelta_{\rB_3}(x) = 
x_1x_2x_3(x_1^2-x_2^2)(x_2^2-x_3^2)(x_3^2-x_1^2).  
\end{equation}
\par
Since $P$ is at least $W(\rA_3)$-symmetric, $\tau_m^{(k)}(x)$ must 
be $W(\rA_3)$-invariant and hence can be written in unique ways as 
polynomials of $e_2(x)$, $e_3(x)$, $e_4(x)$. 
For $m = 1,\dots, 6$, one can write   
\begin{align}
\tau_1^{(k)}(x) &= 0, & \tau_2^{(k)}(x) &= a_2^{(k)} \, e_2(x), \nonumber \\
\tau_3^{(k)}(x) &= a_3^{(k)} \, e_3(x), & 
\tau_4^{(k)}(x) &= a_4^{(k)} \, e_4(x) + b_4^{(k)} \, e_2^2(x), \label{eqn:coeff-tetra} \\ 
\tau_5^{(k)}(x) &= a_5^{(k)} \, e_2(x) \, e_3(x),  & 
\tau_6^{(k)}(x) &= a_6^{(k)} \, e_6(x) + b_6^{(k)} \, e_2(x) \, e_4(x) + c_6^{(k)} \, e_2^3(x).  
\nonumber 
\end{align}
\begin{lemma} \label{lem:pde-tetra}
For $k = 0,1,2,3$, the following hold: 
\begin{enumerate} 
\item If none of $a_2^{(k)}$, $a_3^{(k)}$, $a_4^{(k)}$ is zero, 
then the infinite system \eqref{eqn:pde} is equivalent to the finite 
system \eqref{eqn:pde-a3} so that one has $\H_{P(k)} = \H_{\rA_3}$. 
\item If $a_3^{(k)} = 0$ but none of $a_2^{(k)}$, $a_4^{(k)}$, 
$a_6^{(k)}$ is zero, then the infinite system \eqref{eqn:pde} is 
equivalent to the finite system \eqref{eqn:pde-b3} so that one has  
$\H_{P(k)} = \H_{\rB_3}$. 
\end{enumerate} 
\end{lemma} 
{\it Proof}. First, suppose that none of $a_2^{(k)}$, $a_3^{(k)}$, 
$a_4^{(k)}$ is zero.  
A part of equations \eqref{eqn:coeff-tetra} can then be inverted to 
express $e_2(x)$, $e_3(x)$, $e_4(x)$ as polynomials of 
$\tau_2^{(k)}(x)$, $\tau_3^{(k)}(x)$, $\tau_4^{(k)}(x)$, so that 
system \eqref{eqn:pde-a3} is equivalent to 
$\tau_2^{(k)}(\partial) f = \tau_3^{(k)}(\partial) f = 
\tau_4^{(k)}(\partial) f = 0$.  
For any $m \ge 5$ equation $\tau_m^{(k)}(\partial) f = 0$ is redundant 
because $\tau_m^{(k)}(x)$ is a polynomial of $e_2(x)$, $e_3(x)$, $e_4(x)$. 
This implies that system \eqref{eqn:pde} is equivalent to system 
\eqref{eqn:pde-a3}, leading to the conclusion of assertion (1). 
Next, suppose that $a_3^{(k)} = 0$ but none of $a_2^{(k)}$, 
$a_4^{(k)}$, $a_6^{(k)}$ is zero. 
Another part of equations \eqref{eqn:coeff-tetra} can then be inverted 
to express $e_2(x)$, $e_4(x)$, $e_6(x)$ as polynomials of $\tau_2^{(k)}(x)$, 
$\tau_4^{(k)}(x)$, $\tau_6^{(k)}(x)$, so that system \eqref{eqn:pde-b3} 
is equivalent to $\tau_2^{(k)}(\partial) f = \tau_4^{(k)}(\partial) f = 
\tau_6^{(k)}(\partial) f = 0$. 
Under system \eqref{eqn:pde-b3}, equation $\tau_m^{(k)}(\partial) f = 0$ 
is redundant for $m = 5$ or $m \ge 7$.   
Indeed the $W(\rA_3)$-invariance of $\tau_m^{(k)}(x)$ allows us to write   
\[
\tau_m^{(k)}(\partial) f = \sum_{2 a + 3 b + 4 c = m} \alpha_{abc}^{(k)} 
\cdot e_2^a(\partial) \cdot e_3^b(\partial) \cdot e_4^c(\partial) f,  
\]
with suitable constants $\alpha_{abc}^{(k)}$. 
Equations \eqref{eqn:pde-b3} imply that the summand with index $(a,b,c)$ 
vanishes if either $a \ge 1$, or $b \ge 2$, or $c \ge 1$. 
Thus the index of a nonzero summand, if any, must satisfy $a = c = 0$ and 
$b \le 1$ and so $m = 2 a + 3 b + 4 c \le 2 \cdot 0 + 3 \cdot 1 + 
4 \cdot 0 = 3$. 
Therefore system \eqref{eqn:pde} is equivalent to 
system \eqref{eqn:pde-b3}. 
This proves the second assertion of the lemma.   
\hfill $\Box$
\section{Skeletons of Triakis Tetrahedra} \label{sec:skeleton-tetra} 
Lemma \ref{lem:pde-tetra} is applied to each skeleton of the triakis 
tetrahedron to establish Theorem \ref{thm:main-tetra}. 
\subsection{Vertex Problem} \label{subsec:vertex-tetra}
The polynomial $\tau_m^{(0)}(x)$ in formula \eqref{eqn:tau0} 
is divided into two components: 
\begin{equation} \label{eqn:tau0-tetra}
\tau_m^{(0)}(x) = \tau_{1,m}^{(0)}(x) + \tau_{2,m}^{(0)}(x),
\end{equation}  
according to the two types of vertices, where 
\begin{equation} \label{eqn:tau0-tetra2}
\begin{split}
\tau_{1,m}^{(0)}(x) &= \langle \rA, \, x \rangle^m + \langle \rB, \, x \rangle^m + 
\langle \rC, \, x \rangle^m + \langle \rD, \, x \rangle^m, \\
\tau_{2,m}^{(0)}(x) &= \langle \ra, \, x \rangle^m + \langle \rb, \, x \rangle^m + 
\langle \rc, \, x \rangle^m + \langle \rc, \, x \rangle^m. 
\end{split}
\end{equation} 
Using the coordinates \eqref{eqn:coord-tetra} and the relation 
$\rd = r(\rA+\rB+\rC)/3$ etc., one finds in formulas \eqref{eqn:coeff-tetra},  
\begin{align*}
a_2^{(0)} &= \frac{4}{9} r^2 + 4 > 0, &    
a_3^{(0)} &= \frac{8}{9} (3-r)(r^2 + 3r + 9), \\
a_4^{(0)} &= \frac{16}{81} r^4 + 16 > 0, & 
a_6^{(0)} &= \frac{64}{243}(r^2 + 9)(r^4-9r^2+81) > 0. 
\end{align*}
Note that $a_3^{(0)} = 0$ if and only if $r = 3$. 
Thus if $r = 3$ then assertion (2) of Lemma \ref{lem:pde-tetra} leads to 
the upper case of formula \eqref{eqn:main-tetra} with $k = 0$, while if 
$r \neq 3$ then assertion (1) of that lemma leads to 
the lower case of it.  
In either case there is no gap between $\H_{P(0)}$ and $\H_{G(0)}$.   
\subsection{Edge Problem} \label{subsec:edge-tetra}
It is obvious that for an index $\rA\rB \in I_1$, the foot on the 
line $\ell_{\rA\rB}$ is $\rA\rB = (\rA+\rB)/2$, with this rule 
applying to every index of the same type.  
For indices of the other type in $I_1$, one finds     
\begin{align*}
\rA\rb &= (\b, -\b, -\a), & \rA\rc &= (\b, -\a, -\b), & \rA\rd &= (\a, -\b, -\b),  \\ 
\rB\ra &= (-\b, \, \b, -\a), & \rB\rc &= (-\a, \, \b, -\b), & \rB\rd &= (-\b, \, \a, -\b), \\ 
\rC\ra &= (-\b, -\a, \, \b), & \rC\rb &= (-\a, -\b, \, \b), & \rC\rd &= (-\b, -\b, \, \a), \\
\rD\ra &= (\a, \, \b, \, \b), & \rD\rb &= (\b, \, \a, \, \b), & \rD\rc &= (\b, \, \b, \, \a),  
\end{align*}
where 
\[ 
\a := \frac{4r(r - 3)}{3(r^2 - 2 r + 9)}, \qquad  
\b := \frac{2r(r + 3)}{3(r^2 - 2 r + 9)}. 
\]
For example, formula $\rA\rd = (\a, -\b, -\b)$ follows from the 
condition that the point $\rA\rd$ should lie on the line 
$\ell_{\rA\rd}$, while the vectors $\rA\rd$ and $\rd - \rA$ 
should be orthogonal. 
\par
The polynomial $\tau_m^{(1)}(x)$ in formula \eqref{eqn:tau1} can be 
divided into three components:   
\begin{equation} \label{eqn:tau1-tetra}
\tau_m^{(1)}(x) = \ve_1 \cdot \tau_{1,m}^{(1)}(x) +  \ve_2 \cdot \tau_{2,m}^{(1)}(x) 
+ \ve_3 \cdot \tau_{3,m}^{(1)}(x), 
\end{equation}
according to the three types \eqref{eqn:VE-tetra} of vertex-to-edge 
incidences, where $\tau_{\nu,m}^{(1)}(x)$ is given as in 
Table \ref{tab:tau1-tetra}    
\begin{table}[t]
\[
\begin{split} 
\tau_{1,m}^{(1)}(x) 
&= h_m(\rA, \rA\rB) + h_m(\rA, \rA\rC) + h_m(\rA, \rA\rD) 
+  h_m(\rB, \rA\rB) + h_m(\rB, \rB\rC) + h_m(\rB, \rB\rD) \\
&+ h_m(\rC, \rA\rC) + h_m(\rC, \rB\rC) + h_m(\rC, \rC\rD)
+  h_m(\rD, \rA\rD) + h_m(\rD, \rB\rD) + h_m(\rD, \rC\rD), \\ 
\tau_{2,m}^{(1)}(x) 
&= h_m(\rA, \rA\rb) + h_m(\rA, \rA\rc) + h_m(\rA, \rA\rd) 
+  h_m(\rB, \rB\ra) + h_m(\rB, \rB\rc) + h_m(\rB, \rB\rd) \\ 
&+ h_m(\rC, \rC\ra) + h_m(\rC, \rC\rb) + h_m(\rC, \rC\rd) 
+  h_m(\rD, \rD\ra) + h_m(\rD, \rD\rb) + h_m(\rD, \rD\rc), \\
\tau_{3,m}^{(1)}(x) 
&= h_m(\ra, \rB\ra) + h_m(\ra, \rC\ra) + h_m(\ra, \rD\ra)  
+  h_m(\rb, \rA\rb) + h_m(\rb, \rC\rb) + h_m(\rb, \rD\rb) \\ 
&+ h_m(\rc, \rA\rc) + h_m(\rc, \rB\rc) + h_m(\rc, \rD\rc)  
+  h_m(\rd, \rA\rd) + h_m(\rd, \rB\rd) + h_m(\rd, \rC\rd).  
\end{split}
\]
\vspace{-5mm} 
\caption{The polynomials $\tau_{\nu,m}^{(1)}(x)$, $\nu = 1,2,3$, in 
formula \eqref{eqn:tau1-tetra}.}
\label{tab:tau1-tetra}
\end{table}
and the abbreviation $h_m(\rP, \rQ) := 
h_m(\langle \rP, x \rangle, \, \langle \rQ, x \rangle)$ is used 
for two vectors $\rP$, $\rQ \in \R^3$. 
\par
The three types of vertex-to-edge incidence numbers are 
evaluated as     
\begin{equation} \label{eqn:ve-tetra}
\ve_1 = \sqrt{2}, \qquad 
\ve_2 = \frac{9 - r}{\sqrt{3 (r^2 - 2 r + 9)}}, \qquad
\ve_3 = \frac{r(r - 1)}{\sqrt{3(r^2 - 2 r + 9)}}.  
\end{equation} 
Indeed the formula for $\ve_1$ is easy to see.   
To derive those for $\ve_2$ and $\ve_3$, take a look at the edge 
$\rA\rd$ in Figure \ref{fig:foot}. 
Observing that the unit normal vectors 
$\n_{\rA,\rA\rd}$ and $\n_{\rd,\rA\rd}$ are given by 
$\n_{\rA,\rA\rd} = - \n_{\rd,\rA\rd} = (\rA - \rd)/|\rA - \rd|$, 
where $|\cdot|$ denotes the norm of a vector, one can calculate 
$\ve_2 = [\rA:\rA\rd] = (\rA-\rA\rd)/\n_{\rA,\rA\rd}$ and 
$\ve_3 = [\rd:\rA\rd] = (\rd-\rA\rd)/\n_{\rd,\rA\rd}$ as 
indicated above. 
\par
Putting all these informations together into formulas \eqref{eqn:coeff-tetra}, 
one finds 
\begin{align*}
a_2^{(1)} &= 20 \sqrt{2} + \frac{4}{9}(r^2 + r + 9) 
\sqrt{3 (r^2 - 2r + 9)} > 0, \\
a_3^{(1)} &= 96 \sqrt{2} + \frac{8}{9}(3 - r)(r^2 + 4 r + 9) 
\sqrt{3 (r^2 - 2r + 9)}, \\
a_4^{(1)} &= 48 \sqrt{2} + \frac{16}{81}(r^4 + r^3 - 3 r^2 + 9 r + 81) 
\sqrt{3 (r^2 - 2r + 9)} > 0, \\ 
a_6^{(1)} &= 768 \sqrt{2} + 
\frac{64}{243}(r^6 + r^5 - 3 r^4 - 18 r^3 - 27 r^2 + 81 r + 729) 
\sqrt{3 (r^2 - 2 r + 9)}. 
\end{align*}
Observe that $a_2^{(1)}$ and $a_4^{(1)}$ are positive for every $r > 0$. 
Indeed the former is obvious and the latter follows from the fact that 
$\psi(r) := r^4 + r^3 - 3 r^2 + 9 r + 81$ has a positive value 
$\psi(0) = 81$ at $r = 0$ as well as a positive derivative  
$\psi'(r) = 4 r^3 + 3 (r - 1)^2 + 6$ for every $r > 0$.  
On the other hand, $a_3^{(1)}$ is positive in $0 < r \le 3$,  
strictly decreasing in $r > 3$ and tending to $-\infty$ as $r \to 
+ \infty$. 
Thus there exists a unique positive number $r = r_1$ at which 
$a_3^{(1)} = 0$. 
Observe that 
\[
a_3^{(1)} \left\{\frac{8}{9} (3-r)(r^2 + 4r + 9)
\sqrt{3 (r^2 - 2r + 9)}-96 \sqrt{2}\right\} 
= \frac{64}{27}\, (r^2 -2r+ 3) \, \chi_1(r), 
\]
so $r_1$ must be a positive root of $\chi_1(r)$, where $\chi_1(r)$ is 
the sextic polynomial in \eqref{eqn:sextic}. 
Conversely one can show that $\chi_1(r)$ certainly has a   
unique positive root that yields $r_1 = 3.62398\cdots$.  
Note that $a_6^{(1)} = 1661.36 \cdots$ is nonzero at $r = r_1$.  
Thus if $r = r_1$ then assertion (2) of Lemma \ref{lem:pde-tetra} leads 
to the upper case of formula \eqref{eqn:main-tetra} with $k = 1$, while 
if $r \neq r_1$ then assertion (1) of that lemma leads to the lower case of it. 
There is a gap between $\H_{P(1)}$ and $\H_{G(1)}$ only when $r = r_1$.     
\subsection{Face Problem} \label{subsec:face-tetra}
For each index of $I_2$, the foot on the corresponding affine plane 
is given by  
\begin{align*}
\rA\rB\rc &= (-\d, -\d, -\ga),  &
\rA\rB\rd &= (\d, \, \d, -\ga), &
\rA\rC\rb &= (-\d, -\ga, -\d),  & 
\rA\rC\rd &= (\d, -\ga, \, \d), \\ 
\rA\rD\rb &= (\ga, \, \d, -\d), & 
\rA\rD\rc &= (\ga, -\d, \, \d), & 
\rB\rC\ra &= (-\ga, -\d, -\d),  &
\rB\rC\rd &= (-\ga, \, \d, \, \d), \\
\rB\rD\ra &= (\d, \, \ga, -\d), &
\rB\rD\rc &= (-\d, \, \ga, \d), &  
\rC\rD\ra &= (\d, -\d, \, \ga), & 
\rC\rD\rb &= (-\d, \, \d, \ga),  
\end{align*}
where
\[
\gamma := \frac{2 r^2}{3(r^2 - 2r + 3)}, \qquad 
\delta := \frac{r(r - 3)}{3(r^2 - 2r + 3)}. 
\]
For example, formula $\rA\rB\rd = (\d, \d, -\ga)$ follows from the 
condition that the point $\rA\rB\rd$ should lie on the plane  
$H_{\rA\rB\rd}$, while the vector $\rA\rB\rd$ should be 
orthogonal to both $\rd - \rA$ and $\rd - \rB$. 
\par
Observe that up to symmetries there are three types of 
vertex-edge-face flags: 
\[
\begin{matrix}
(1) \quad \rA \prec \rA\rB \prec \rA\rB\rd \quad (\VE 1) \,\, \& \,\, (\EF 1), \quad   & 
(2) \quad \rA \prec \rA\rd \prec \rA\rB\rd \quad (\VE 2) \,\, \& \,\, (\EF 2), \\[1mm] & 
(3) \quad \rd \prec \rA\rd \prec \rA\rB\rd \quad (\VE 3) \,\, \& \,\, (\EF 2), 
\end{matrix}
\]
according to which the polynomial $\tau_m^{(2)}(x)$ in formula \eqref{eqn:tau2} 
can be divided into three components:    
\begin{equation} \label{eqn:tau2-tetra}
\tau_m^{(2)}(x) = \ve_1 \cdot \ef_1 \cdot \tau_{1,m}^{(2)}(x) + 
\ve_2 \cdot \ef_2 \cdot \tau_{2,m}^{(2)}(x) + 
\ve_3 \cdot \ef_2 \cdot \tau_{3,m}^{(2)}(x), 
\end{equation}
where the polynomials $\tau_{\nu,m}^{(2)}(x)$, $\nu = 1,2,3$, are given as in 
Table \ref{tab:tau2-tetra}     
\begin{table}[t]
\[
\begin{split}
\tau_{1,m}^{(2)}(x) 
&= h_m(\rA, \rA\rB, \rA\rB\rc) + h_m(\rA, \rA\rB, \rA\rB\rd) + h_m(\rA, \rA\rC, \rA\rC\rb) + h_m(\rA, \rA\rC, \rA\rC\rd) \\
&+ h_m(\rA, \rA\rD, \rA\rD\rb) + h_m(\rA, \rA\rD, \rA\rD\rc) + h_m(\rB, \rA\rB, \rA\rB\rc) + h_m(\rB, \rA\rB, \rA\rB\rd) \\
&+ h_m(\rB, \rB\rC, \rB\rC\ra) + h_m(\rB, \rB\rC, \rB\rC\rd) + h_m(\rB, \rB\rD, \rB\rD\ra) + h_m(\rB, \rB\rD, \rB\rD\rc) \\
&+ h_m(\rC, \rA\rC, \rA\rC\rb) + h_m(\rC, \rA\rC, \rA\rC\rd) + h_m(\rC, \rB\rC, \rB\rC\ra) + h_m(\rC, \rB\rC, \rB\rC\rd) \\
&+ h_m(\rC, \rC\rD, \rC\rD\ra) + h_m(\rC, \rC\rD, \rC\rD\rb) + h_m(\rD, \rA\rD, \rA\rD\rb) + h_m(\rD, \rA\rD, \rA\rD\rc) \\
&+ h_m(\rD, \rB\rD, \rB\rD\ra) + h_m(\rD, \rB\rD, \rB\rD\rc) + h_m(\rD, \rC\rD, \rC\rD\ra) + h_m(\rD, \rC\rD, \rC\rD\rb), \\
\tau_{2,m}^{(2)}(x) 
&= h_m(\rA, \rA\rb, \rA\rC\rb) + h_m(\rA, \rA\rb, \rA\rD\rb) + h_m(\rA, \rA\rc, \rA\rB\rc) + h_m(\rA, \rA\rc, \rA\rD\rc) \\
&+ h_m(\rA, \rA\rd, \rA\rB\rd) + h_m(\rA, \rA\rd, \rA\rC\rd) + h_m(\rB, \rB\ra, \rB\rC\ra) + h_m(\rB, \rB\ra, \rB\rD\ra) \\
&+ h_m(\rB, \rB\rc, \rA\rB\rc) + h_m(\rB, \rB\rc, \rB\rD\rc) + h_m(\rB, \rB\rd, \rA\rB\rd) + h_m(\rB, \rB\rd, \rB\rC\rd) \\
&+ h_m(\rC, \rC\ra, \rB\rC\ra) + h_m(\rC, \rC\ra, \rC\rD\ra) + h_m(\rC, \rC\rb, \rA\rC\rb) + h_m(\rC, \rC\rb, \rC\rD\rb) \\
&+ h_m(\rC, \rC\rd, \rA\rC\rd) + h_m(\rC, \rC\rd, \rB\rC\rd) + h_m(\rD, \rD\ra, \rB\rD\ra) + h_m(\rD, \rD\ra, \rC\rD\ra) \\
&+ h_m(\rD, \rD\rb, \rA\rD\rb) + h_m(\rD, \rD\rb, \rC\rD\rb) + h_m(\rD, \rD\rc, \rA\rD\rc) + h_m(\rD, \rD\rc, \rB\rD\rc), \\ 
\tau_{3,m}^{(2)}(x)
&= h_m(\ra, \rB\ra, \rB\rC\ra) + h_m(\ra, \rB\ra, \rB\rD\ra) + h_m(\ra, \rC\ra, \rB\rC\ra) + h_m(\ra, \rC\ra, \rC\rD\ra) \\
&+ h_m(\ra, \rD\ra, \rB\rD\ra) + h_m(\ra, \rD\ra, \rC\rD\ra) + h_m(\rb, \rA\rb, \rA\rC\rb) + h_m(\rb, \rA\rb, \rA\rD\rb) \\
&+ h_m(\rb, \rC\rb, \rA\rC\rb) + h_m(\rb, \rC\rb, \rC\rD\rb) + h_m(\rb, \rD\rb, \rA\rD\rb) + h_m(\rb, \rD\rb, \rC\rD\rb) \\
&+ h_m(\rc, \rA\rc, \rA\rB\rc) + h_m(\rc, \rA\rc, \rA\rD\rc) + h_m(\rc, \rB\rc, \rA\rB\rc) + h_m(\rc, \rB\rc, \rB\rD\rc) \\
&+ h_m(\rc, \rD\rc, \rA\rD\rc) + h_m(\rc, \rD\rc, \rB\rD\rc) + h_m(\rd, \rA\rd, \rA\rB\rd) + h_m(\rd, \rA\rd, \rA\rC\rd) \\
&+ h_m(\rd, \rB\rd, \rA\rB\rd) + h_m(\rd, \rB\rd, \rB\rC\rd) + h_m(\rd, \rC\rd, \rA\rC\rd) + h_m(\rd, \rC\rd, \rB\rC\rd).   
\end{split} 
\]
\vspace{-5mm}
\caption{The polynomials $\tau_{\nu,m}^{(2)}(x)$, $\nu = 1,2,3$, in 
formula \eqref{eqn:tau2-tetra}.}
\label{tab:tau2-tetra} 
\end{table} 
and the abbreviation $h_m(\rP,\rQ,\rR) := 
h_m(\langle \rP, x \rangle, \, \langle \rQ, x \rangle, \, \langle \rR, x \rangle)$ is 
used for three vectors $\rP$, $\rQ$, $\rR \in \R^3$. 
\par
While the vertex-to-edge incidence numbers are given in 
\eqref{eqn:ve-tetra}, the edge-to-face ones are  
\begin{equation*} \label{eqn:ef-tetra}
\ef_1 = \frac{3 - r}{\sqrt{3 (r^2 - 2r + 3)}}, \qquad 
\ef_2 = \frac{\sqrt{2} \, r(r - 1)}{\sqrt{(r^2 - 2r + 3) (r^2 - 2r + 9)}}.     
\end{equation*}
So upon multiplying by a nonzero constant simultaneously, one may put 
\begin{align*}
\ve_1 \cdot \ef_1 &= \frac{3 - r}{r^2- 2r + 3}, \qquad &   
\ve_2 \cdot \ef_2 &= \frac{(9 - r) r (r - 1)}{(r^2 -2r + 3)(r^2 -2r + 9)}, \\  
                  &                                   &
\ve_3 \cdot \ef_2 &= \frac{r^2(r - 1)^2}{(r^2 -2r + 3)(r^2 -2r + 9)}.   
\end{align*}  
Notice that what is important in expression \eqref{eqn:tau2-tetra} is only 
the ratio $\ve_1 \cdot \ef_1 : \ve_2 \cdot \ef_2 : \ve_3 \cdot \ef_2$. 
\par
Putting all these informations together into formulas \eqref{eqn:coeff-tetra}, 
one finds  
\begin{align*}
a_2^{(2)} &= \frac{8}{3} (r^2 + 2r + 15) > 0, &  
a_3^{(2)} &= \frac{16}{3} (3 - r)  (r^2 + 5r + 12), \\ 
a_4^{(2)} &= \frac{32}{27} (r^4 + 2 r^3 - 3 r^2 + 81) > 0, & 
a_6^{(2)} &= \frac{128}{81} (r^6+ 2 r^5- 3 r^4- 27 r^3- 54 r^2 + 81 r + 972).  
\end{align*}
It is easy to see that $a_2^{(2)}$ and $a_4^{(2)}$ are positive 
for every $r \in \R$ and that $a_3^{(2)} = 0$ if and only if 
$r = 3$, in which case $a_6^{(2)} = 1536$ is nonzero. 
Thus if $r = 3$ then assertion (2) of Lemma \ref{lem:pde-tetra} yields 
the upper case of formula \eqref{eqn:main-tetra} with $k = 2$, while if 
$r \neq 3$ then assertion (1) of that lemma yields the lower case of it.   
In either case there is no gap between $\H_{P(2)}$ and $\H_{G(2)}$. 
\par
By \cite[Theorem 2.2]{Iwasaki1} the volume problem ($k = 3$) has the same 
solution as the face problem ($k = 2$), since $P$ is isohedral. 
The proof of Theorem \ref{thm:main-tetra} is now complete.    
\section{A Family of Isohedral Triakis Octahedra} \label{sec:family-octa}
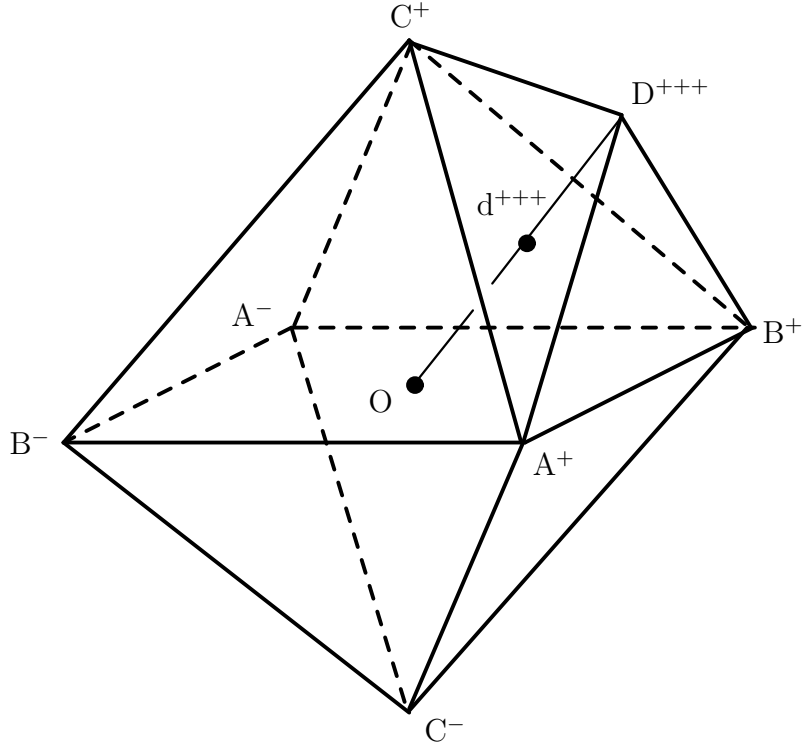
\begin{figure}[t]
\begin{center}
\unitlength 0.1in
\begin{picture}( 39.0000, 37.5000)( 17.3000,-40.1000)
%
\special{pn 20}%
\special{pa 2010 2600}%
\special{pa 4400 2600}%
\special{fp}%
%
\special{pn 20}%
\special{pa 3200 2000}%
\special{pa 5590 2000}%
\special{da 0.070}%
%
\special{pn 20}%
\special{pa 3190 2000}%
\special{pa 2010 2600}%
\special{da 0.070}%
%
\special{pn 20}%
\special{pa 5580 2000}%
\special{pa 4400 2600}%
\special{fp}%
%
\special{pn 20}%
\special{pa 3800 4000}%
\special{pa 5550 2000}%
\special{fp}%
%
\special{pn 20}%
\special{pa 4390 2600}%
\special{pa 3800 3990}%
\special{fp}%
%
\special{pn 20}%
\special{pa 2010 2600}%
\special{pa 3800 4010}%
\special{fp}%
%
\special{pn 20}%
\special{pa 3190 2000}%
\special{pa 3790 4000}%
\special{da 0.070}%
%
\special{pn 20}%
\special{pa 3800 500}%
\special{pa 4380 2600}%
\special{fp}%
%
\special{pn 20}%
\special{pa 3800 500}%
\special{pa 5560 2000}%
\special{da 0.070}%
%
\special{pn 20}%
\special{pa 3800 500}%
\special{pa 2010 2600}%
\special{fp}%
%
\special{pn 20}%
\special{pa 3810 520}%
\special{pa 3200 2000}%
\special{da 0.070}%
%
\special{pn 20}%
\special{pa 4900 890}%
\special{pa 4390 2600}%
\special{fp}%
%
\special{pn 20}%
\special{pa 3800 510}%
\special{pa 4900 890}%
\special{fp}%
%
\special{pn 20}%
\special{pa 4900 900}%
\special{pa 5570 2000}%
\special{fp}%
%
\special{pn 13}%
\special{pa 4900 900}%
\special{pa 4230 1770}%
\special{fp}%
%
\special{pn 13}%
\special{pa 4130 1910}%
\special{pa 3830 2290}%
\special{fp}%
%
\special{pn 20}%
\special{sh 1.000}%
\special{ar 4410 1560 36 36  0.0000000 6.2831853}%
%
\special{pn 20}%
\special{sh 1.000}%
\special{ar 3830 2300 36 36  0.0000000 6.2831853}%
\put(44.4000,-27.7000){\makebox(0,0)[lb]{$\rA^+$}}%
\put(56.3000,-20.8000){\makebox(0,0)[lb]{$\rB^+$}}%
\put(37.0000,-4.3000){\makebox(0,0)[lb]{$\rC^+$}}%
\put(49.5000,-8.3000){\makebox(0,0)[lb]{$\rD^{+++}$}}%
\put(35.9000,-24.4000){\makebox(0,0)[lb]{$\rO$}}%
\put(38.8000,-41.6000){\makebox(0,0)[lb]{$\rC^-$}}%
\put(17.3000,-26.6000){\makebox(0,0)[lb]{$\rB^-$}}%
\put(28.8000,-19.9000){\makebox(0,0)[lb]{$\rA^-$}}%
\put(41.5000,-14.1000){\makebox(0,0)[lb]{$\rd^{+++}$}}%
\end{picture}%
\end{center}
\caption{Adjoining a pyramid to a face of a regular octahedron; 
$\overline{\rO\rD^{abc}} : \overline{\rO\rd^{abc}} = r : 1$.}
\label{fig:triakis-octa}
\end{figure}

Let $\Omega$ be a regular octahedron with center at the origin $\rO$. 
The six vertices of $\Omega$ can be written $\rA^{\pm}$, $\rB^{\pm}$, 
$\rC^{\pm}$, where $\rA^+$ and $\rA^-$ are antipodal to each other 
and so on. 
The eight faces are then given by $\rA^a\rB^b\rC^c$ with $a, b, c = \pm$. 
An {\sl isohedral triakis octahedron} $P$ is obtained from $\Omega$ by 
adjoining to each face of it a pyramid of appropriate height, or 
excavating such a pyramid, just as in the construction of an isohedral 
triakis tetrahedron. 
Let $\rD^{abc}$ denote the top vertex of the pyramid based on the 
face $\rA^a\rB^b\rC^c$. 
The polyhedron $P$ depends on a parameter $r > 0$ in such a manner 
that the distance ratio $\overline{\rO\rD^{abc}} : \overline{\rO\rd^{abc}}$ 
is $r : 1$, where $\rd^{abc}$ is the center of the face $\rA^a\rB^b\rC^c$ 
(see Figure \ref{fig:triakis-octa}).  
\par
\begin{figure}[t]
\begin{minipage}{0.49\hsize}
\begin{center}
\includegraphics*[width=60mm,height=60mm]{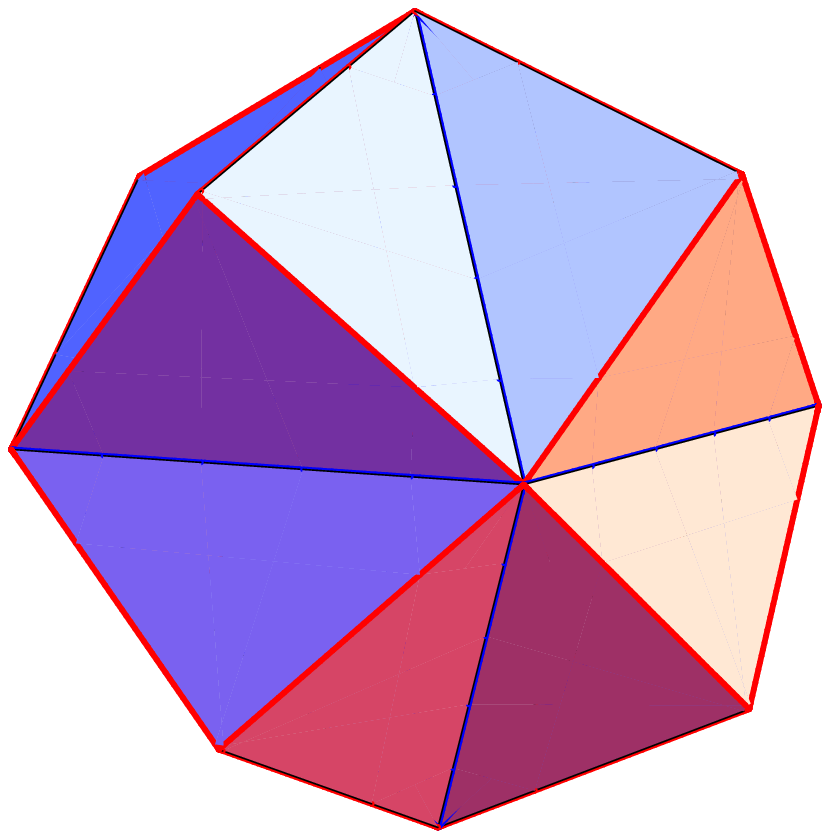} 
\end{center}
\end{minipage}
\hspace{-5mm}
\begin{minipage}{0.49\hsize}
\begin{center}
\includegraphics*[width=60mm,height=60mm]{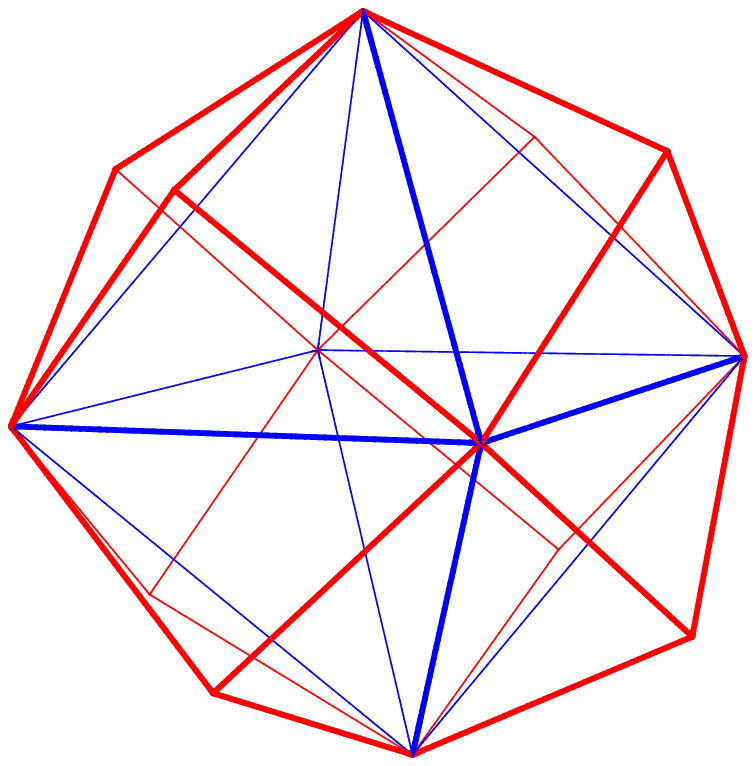} 
\end{center}
\end{minipage}
\caption{The isohedral triakis octahedron with $r = 1.82977\cdots$.}
\label{fig:triakis-octa2}
\end{figure}

Let us look more closely at the polyhedron $P$ for various values 
of $r$. 
The values $r = 1$ and $r = 3/2$ are special in that as a point set,  
$P$ degenerates to the original octahedron $\Omega$ at $r = 1$ and to 
a rhombic dodecahedron at $r = 3/2$ where two neighboring faces, say, 
$\rA^+\rB^+\rD^{+++}$ and $\rA^+\rB^+\rD^{++-}$ are coplanar.  
Note that $P$ is convex if and only if $1 \le r \le 3/2$, in which 
interval the value $r = 3(\sqrt{2}-1) = 1.24264\cdots$ is distinguished 
in that $P$ becomes a Catalan solid \cite{Catalan,KKK}, or an Archimedean dual 
solid, called the small triakis octahedron, whose dual is the truncated 
cube $[3,8,8]$. 
As in the case of triakis tetrahedron we employ the convention that 
as a combinatorial polyhedron, $P$ has the constant skeletal 
structure for all $r > 0$, even at $r = 1$ and $r = 3/2$. 
Figure \ref{fig:triakis-octa2} exhibits the shape of $P$ when  
$r = 1.82977\cdots$; see the $k = 2, 3$ case of 
Theorem \ref{thm:main-octa} for the origin of this particular value.      
\par
If the vertices of our octahedron $\Omega$ are taken as 
\begin{equation} \label{eqn:coord-octa}
\rA^{\pm} = (\pm1, 0, 0), \qquad \rB^{\pm} = (0, \pm1, 0), \qquad 
\rC^{\pm} = (0, 0, \pm1),  
\end{equation}
then the symmetry group of $\Omega$ is the same as that of the cube 
$C$ in \S\ref{sec:family-tetra}, namely, the Weyl group $W(\rB_3)$. 
It is just the symmetry group $G = G(k)$ of $P(k)$ for every 
$k = 0, 1, 2, 3$ and $r > 0$.   
\begin{theorem} \label{thm:main-octa} 
For the family of isohedral triakis octahedra, $\H_{G(k)} \varsubsetneq 
\H_{P(k)}$ if and only if 
\begin{itemize}
\item $k = 0$ and $r = r_0 := 3 \cdot 2^{-3/4} = 1.78381\cdots$,    
\item $k = 1$ and $r = r_1 := 2.24580\cdots$ is the unique positive root 
of an octic equation 
\begin{equation} \label{eqn:octic}
\chi_1(r) := 16 r^8 + 32 r^7 + 40 r^6 - 48 r^5  - 396 r^4 - 432 r^3  - 810 r^2  - 972 r - 729 = 0,  
\end{equation}
\item $k = 2, 3$ and $r = r_2 = r_3 := 1.82977\cdots$ is the unique 
positive root of a quartic equation 
\begin{equation} \label{eqn:quartic}
\chi_2(r) := 4 r^4  + 8 r^3  + 6 r^2  - 18 r - 81 = 0.  
\end{equation}
\end{itemize}
If any of these is the case, then $\dim \H_{G(k)} = 48 < \dim \H_{P(k)} = 96$ 
and as $\R[\partial]$-modules, $\H_{G(k)} = \H_{\rB_3}$ is generated 
by the polynomial $\varDelta_{\rB_3}(x)$ in formula \eqref{eqn:fap-b3} 
while $\H_{P(k)}$ is generated by 
\begin{equation} \label{eqn:F}
F(x) := \varDelta_{\rB_3}(x) \cdot 
\{5(x_1^4+x_2^4+x_3^4)-13(x_1^2 x_2^2 + x_2^2x_3^2+x_3^2x_1^2) \}. 
\end{equation}
\end{theorem}
\par
As in the case of Theorem \ref{thm:main-tetra} this theorem is also 
established through the analysis of PDEs \eqref{eqn:pde}. 
To adjust the notation to the current situation we represent the index 
sets $I_0$, $I_1$, $I_2$ in \S\ref{sec:inv-pde} by letting the vertices, 
edges and faces of $P$ to speak of themselves, that is, by setting    
\begin{align*}
I_0 &= \{ \rA^a, \, \rB^b, \, \rC^c \} \cup \{ \rD^{abc} \},  \\
I_1 &= \{\rA^a\rB^b, \, \rA^a\rC^c, \, \rB^b\rC^c \} \cup 
\{\rA^a\rD^{abc}, \, \rB^b\rD^{abc}, \, \rC^c\rD^{abc} \}, \\
I_2 &= \{\rA^a\rB^b\rD^{abc}, \, \rA^a\rC^c\rD^{abc}, \, \rB^b\rC^c\rD^{abc} \},   
\end{align*}
where $a$, $b$, $c$ take $\pm$ signs freely and $\{ \cdots \}$ stands 
for an orbit under symmetries.   
For an index $\rA^a\rB^b \in I_1$ the same symbol $\rA^a\rB^b$ denotes 
the foot of orthogonal projection from the origin $\rO$ to the affine 
line passing through $\rA^a$ and $\rB^b$; this rule also 
applies to another index $\rA^a\rD^{abc} \in I_1$ as well as to all the 
other indices of $I_1$. 
Similarly, for an index $\rA^a\rB^b\rD^{abc} \in I_2$ the 
same symbol $\rA^a\rB^b\rD^{abc}$ denotes the foot of orthogonal 
projection from $\rO$ to the affine plane passing through 
$\rA^a$, $\rB^b$, $\rD^{abc}$, with this rule applying to all the 
other indices of $I_2$; see Figure \ref{fig:triakis-octa}.    
\par
Up to symmetries there are three types of vertex-to-edge 
incidence relations: 
\[
(\VE 1) \quad \rA^+ \prec \rA^+\rB^+, \qquad 
(\VE 2) \quad \rA^+ \prec \rA^+\rD^{+++}, \qquad 
(\VE 3) \quad \rD^{+++} \prec \rA^+\rD^{+++}.  
\]
On the other hand, up to symmetries there are two types of 
edge-to-face incidence relations:  
\[
(\EF 1) \quad \rA^+\rB^+ \prec \rA^+\rB^+\rD^{+++}, \qquad 
(\EF 2) \quad \rA^+\rD^{+++} \prec \rA^+\rB^+\rD^{+++}.  
\]
The incidence number of an incidence relation depends only on 
its type. 
Let $\ve_{\nu}$ denote the incidence number of type ($\VE \nu$), 
$\nu = 1,2,3$, and $\ef_{\nu}$ denote that of type ($\EF \nu$), 
$\nu = 1,2$, respectively.   
Notice that the presentation around here is quite similar to the 
presentation around formulas \eqref{eqn:VE-tetra} and \eqref{eqn:EF-tetra}. 
Indeed the treatment of triakis octahedra will be very parallel  
to that of triakis tetrahedra, so that the main focus in what 
follows will be on how one can modify the arguments in the 
tetrahedral case to the octahedral one.   
\par
We consider how Lemma \ref{lem:pde-tetra} can be modified. 
Since $P$ is $W(\rB_3)$-symmetric, $\tau_m^{(k)}(x)$ must be 
$W(\rB_3)$-invariant and hence can be written in unique ways as 
polynomials of $e_2(x)$, $e_4(x)$, $e_6(x)$. 
Note that $\tau_m^{(k)}(x)$ is identically zero whenever $m$ is odd.  
For $m = 2, 4, 6, 8$, one can write   
\begin{equation} \label{eqn:coeff-octa}
\begin{split}
\tau_2^{(k)}(x) &= a_2^{(k)} \, e_2(x), \\ 
\tau_4^{(k)}(x) &= a_4^{(k)} \, e_4(x) + b_4^{(k)} \, e_2^2(x), \\ 
\tau_6^{(k)}(x) &= a_6^{(k)} \, e_6(x) + b_6^{(k)} \, e_2(x) \, e_4(x) + 
c_6^{(k)} \, e_2^3(x), \\
\tau_8^{(k)}(x) &= a_8^{(k)} \, e_4^2(x) + b_8^{(k)} \, e_2(x) \, e_6(x) + 
c_8^{(k)} \, e_2^2(x) \, e_4(x) + d_8^{(k)} \, e_2^4(x). 
\end{split}
\end{equation}
\begin{lemma} \label{lem:pde-octa}
For $k = 0, 1, 2, 3$, the following hold: 
\begin{enumerate} 
\item If none of $a_2^{(k)}$, $a_4^{(k)}$, $a_6^{(k)}$ is zero, 
then the infinite system \eqref{eqn:pde} is equivalent to the finite 
system \eqref{eqn:pde-b3} so that one has $\H_{P(k)} = \H_{\rB_3}$. 
\item If $a_4^{(k)} = 0$ but none of $a_2^{(k)}$, $a_6^{(k)}$, 
$a_8^{(k)}$ is zero, then system \eqref{eqn:pde} is equivalent to 
\begin{equation} \label{eqn:pde-octa}
e_2(\partial) f = e_6(\partial) f = e_4^2(\partial) f = 0. 
\end{equation} 
If $\Sol$ denotes the solution space to system \eqref{eqn:pde-octa}, 
then there exists an exact sequence 
\begin{equation} \label{eqn:exact}
\begin{CD}
0 @>>> \H_{\rB_3} @> {\scriptstyle \mathrm{inclusion}} >> 
\Sol @> e_4(\partial) >> \H_{\rB_3} @>>> 0. 
\end{CD} 
\end{equation}
As an $\R[\partial]$-module, $\Sol$ is generated by the polynomial 
$F(x)$ defined in formula \eqref{eqn:F}. 
\end{enumerate} 
\end{lemma}
{\it Proof}. The proof of assertion (1) of this lemma is the same 
as that of assertion (2) of Lemma \ref{lem:pde-tetra}. 
Under the assumption of assertion (2) of this lemma, the first, third 
and fourth equations of \eqref{eqn:coeff-octa} imply that 
system \eqref{eqn:pde-octa} leads to $\tau_2^{(k)}(\partial) f = 
\tau_6^{(k)}(\partial) f = \tau_8^{(k)}(\partial) f = 0$ and 
conversely the latter leads back to the former. 
Under system \eqref{eqn:pde-octa}, equation $\tau_m^{(k)}(\partial) f = 0$ 
is redundant for $m = 4$ or for any even $m \ge 10$. 
Indeed, for $m = 4$ it is trivial from $a_4^{(k)} = 0$ and the second 
equation of \eqref{eqn:coeff-octa}. 
For $m \ge 10$ the $W(\rB_3)$-invariance of $\tau_m^{(k)}(x)$ allows us 
to write   
\[
\tau_m^{(k)}(\partial) f = \sum_{2 a + 4 b + 6 c = m} \alpha_{abc}^{(k)} 
\cdot e_2^a(\partial) \cdot e_4^b(\partial) \cdot e_6^c(\partial) f,  
\]
with suitable constants $\alpha_{abc}^{(k)}$. 
Equations \eqref{eqn:pde-octa} imply that the summand with index $(a,b,c)$ 
vanishes if either $a \ge 1$, or $b \ge 2$, or $c \ge 1$. 
Thus the index of a nonzero summand, if any, must satisfy $a = c = 0$ and 
$b \le 1$ and so $m = 2 a + 4 b + 6 c \le 2 \cdot 0 + 4 \cdot 1 + 
4 \cdot 0 = 4$.  
Therefore system \eqref{eqn:pde} is equivalent to system \eqref{eqn:pde-octa}. 
Next we verify the exact sequence \eqref{eqn:exact}. 
The inclusion $\H_{\rB_3} \subset \Sol$ and the well-definedness of 
$e_4(\partial) : \Sol \to \H_{\rB_3}$ are obvious since $\H_{\rB_3}$  
is the solution space to system \eqref{eqn:pde-b3} while $\Sol$ is to 
system \eqref{eqn:pde-octa}.  
For the same reason sequence \eqref{eqn:exact} is exact at the middle 
term $\Sol$.  
A direct check shows that polynomial $F(x)$ in formula \eqref{eqn:F} 
satisfies 
\[
e_2(\partial) F = e_6(\partial) F = 0, \qquad e_4(\partial) F = - 15120 
\varDelta_{\rB_3}. 
\]
Hence $F \in \Sol$ and $e_4(\partial)$ sends $F$ to $\varDelta_{\rB_3}$ 
up to a nonzero constant multiple. 
Thus $e_4(\partial) : \Sol \to \H_{\rB_3}$ is surjective, since 
it is an $\R[\partial]$-homomorphism and $\varDelta_{\rB_3}$ 
generates $\H_{\rB_3}$ as an $\R[\partial]$-module. 
Finally we show that $F$ generates $\Sol$ as an $\R[\partial]$-module. 
For any $f \in \Sol$, consider $e_4(\partial) f \in \H_{\rB_3}$. 
There is a polynomial $\varphi(x)$ such that $e_4(\partial) f 
= -15120 \, \varphi(\partial) \varDelta_{\rB_3}$. 
Then $g := f - \varphi(\partial) F \in \Sol$ satisfies 
$e_4(\partial) g = e_4(\partial) f - \varphi(\partial) e_4(\partial) F 
= e_4(\partial) f + 15120 \, \varphi(\partial) \varDelta_{\rB_3} = 0$. 
Exact sequence \eqref{eqn:exact} tells us that $g \in \H_{\rB_3}$ and 
so there is a polynomial $\psi(x)$ such that 
$g = -15120 \, \psi(\partial) \varDelta_{\rB_3} = \psi(\partial) 
e_4(\partial) F$. 
Now if $\eta(x) := \varphi(x) + \psi(x) e_4(x)$ then 
$f = \eta(\partial) F$.  
This proves the last claim. \hfill $\Box$ 
\section{Skeletons of Triakis Octahedra} \label{sec:skeleton-octa} 
Lemma \ref{lem:pde-octa} is applied to each skeleton of the triakis 
octahedron to establish Theorem \ref{thm:main-octa}.  
\subsection{Vertex Problem} \label{subsec:vertex-octa}
Formula \eqref{eqn:tau0-tetra} for triakis tetrahedra carries 
over to the case of triakis octahedra, but this time formula 
\eqref{eqn:tau0-tetra2} should be replaced by 
\[
\tau_{1,m}^{(0)}(x) = \sum_a \langle \rA^a, \, x \rangle^m + 
\sum_b \langle \rB^b, \, x \rangle^m + \sum_c \langle \rC^c, \, x \rangle^m, 
\quad 
\tau_{2,m}^{(0)}(x) = \sum_{(a,b,c)} \langle \rD^{abc}, \, x \rangle^m,  
\]
where the sums are taken over all $a, b, c = \pm$. 
Using the coordinates \eqref{eqn:coord-octa} and the relation 
$\rD^{abc} = r(\rA^a + \rB^b + \rC^c)/3$, one finds in formulas  
\eqref{eqn:coeff-octa},    
\begin{align*}
a_2^{(0)} &= \frac{8}{9} r^2 + 2 > 0, &    
a_4^{(0)} &= \frac{32}{81} \left( r^4 - \frac{3^4}{2^3} \right), \\
a_6^{(0)} &= \frac{2}{243}(4 r^2 + 9)(16 r^4 -36 r^2 + 81) > 0, & 
a_8^{(0)} &= \frac{128}{6561} r^8 + 4 > 0. 
\end{align*} 
Note that $a_4^{(0)} = 0$ if and only if $r = r_0 := 3 \cdot 2^{-3/4} 
= 1.78381\cdots$. 
Thus if $r = r_0$ then assertion (2) of Lemma \ref{lem:pde-octa} 
yields $\H_{P(0)} = \Sol$, while if $r \neq r_0$ then assertion (1) 
of Lemma \ref{lem:pde-octa} yields $\H_{P(0)} = \H_{\rB_3}$.    
One has $\H_{G(0)} \varsubsetneq \H_{P(0)}$ only when $r = r_0$. 
This proves Theorem \ref{thm:main-octa} for $k = 0$.    
\subsection{Edge Problem} \label{subsec:edge-octa}
It is obvious that for an index $\rA^a\rB^b \in I_1$, the corresponding 
foot is $\rA^a\rB^b = (\rA^a+\rB^b)/2$ with this rule applying to every 
index of the same type.  
For indices of the other type in $I_1$, one finds     
\[
\rA^a\rD^{abc} = (a u, \, b v, \, c v), \qquad  
\rB^b\rD^{abc} = (a v, \, b u, \, c v), \qquad  
\rC^c\rD^{abc} = (a v, \, b v, \, c u),   
\]
with $a, b, c = \pm$, where 
\[ 
u := \frac{2 r^2}{3(r^2 - 2 r + 3)}, \qquad  
v := \frac{r (3 - r)}{3 (r^2 - 2 r + 3)}. 
\]
\par
Formula \eqref{eqn:tau1-tetra} remains true for triakis octahedra if 
formulas in Table \ref{tab:tau1-tetra} are replaced by 
\[
\begin{split}
\tau_{1,m}^{(1)}(x)  
&= \sum h_m(\rA^a, \rA^a\rB^b) + \sum h_m(\rA^a, \rA^a\rC^c) + \sum h_m(\rB^b, \rA^a\rB^b) \\ 
&+ \sum h_m(\rB^b, \rB^b\rC^c) + \sum h_m(\rC^c, \rA^a\rC^c) + \sum h_m(\rC^c, \rB^b\rC^c), \\
\tau_{2,m}^{(1)}(x)  
&= \sum h_m(\rA^a, \rA^a\rD^{abc}) + \sum h_m(\rB^b, \rB^b\rD^{abc}) + 
\sum h_m(\rC^c, \rC^c\rD^{abc}), \\
\tau_{3,m}^{(1)}(x)  
&= \sum h_m(\rD^{abc}, \rA^a\rD^{abc}) + \sum h_m(\rD^{abc}, \rB^b\rD^{abc}) + 
\sum h_m(\rD^{abc}, \rC^c\rD^{abc}).      
\end{split} 
\] 
The three types of vertex-to-edge incidence numbers are 
evaluated as     
\begin{equation} \label{eqn:ve-octa}
\ve_1 = \frac{1}{\sqrt{2}},  
\qquad 
\ve_2 = \frac{3 - r}{\sqrt{3 (r^2 - 2 r + 3)}}, 
\qquad
\ve_3 = \frac{r(r - 1)}{\sqrt{3 (r^2 - 2 r + 3)}}.   
\end{equation} 
\par
Putting all these informations together into formulas \eqref{eqn:coeff-octa}, 
one finds 
\begin{align*}
a_2^{(1)} &= 8 \sqrt{2} +\frac{8}{9}(r^2+r+3) \sqrt{3(r^2-2r+3)} > 0, \\    
a_4^{(1)} &= -12 \sqrt{2} + \frac{16}{81} (2r^4+2r^3-9r-27) \sqrt{3(r^2-2r+3)}, \\
a_6^{(1)} &= 12 \sqrt{2} + \frac{8}{243} (16r^6+16r^5-18r^3+81r+243) \sqrt{3(r^2-2r+3)} > 0, \\ 
a_8^{(1)} &= 12 \sqrt{2} + \frac{16}{6561}(8r^8+8r^7-36r^5-108r^4-162r^3+729r+2187) \sqrt{3(r^2-2r+3)}. 
\end{align*}
Observe that $a_2^{(1)}$ and $a_6^{(1)}$ are positive for every $r > 0$. 
Indeed the former is obvious and the latter follows from the fact that 
$\psi(r) := 16r^6+16r^5-18r^3+81r+243$ has a positive value 
$\psi(0) = 243$ at $r = 0$ and a positive derivative  
$\psi'(r) = 96 r^5 + 71 r^4 + (3 r^2 - 9)^2$ for every $r > 0$.  
On the other hand, there exists a unique positive number $r = r_1$ at which 
$a_4^{(1)} = 0$, because 
\[
a_4^{(1)} = -4 (4 + 3 \sqrt{2}) < 0 \quad \mbox{at \,\, $r = 0$}; \qquad 
\frac{d a_4^{(1)}}{dr} = \frac{160 r^3 (r^2 - r + 1)}{27 \sqrt{3(r^2 - 2 r + 3)}} > 0 
\quad \mbox{for \,\, $r > 0$},   
\]
and $a_4^{(1)}$ tends to $+\infty$ as $r \to + \infty$. 
If $\chi_1(r)$ is the octic polynomial defined in formula
\eqref{eqn:octic},   
\[
a_4^{(1)}\left\{ 12 \sqrt{2} + \frac{16}{81} (2r^4+2r^3-9r-27) \sqrt{3(r^2-2r+3)} \right\} 
= \frac{32}{2187}(2r^2-4r+3) \, \chi_1(r).     
\]
So $r_1 = 2.24580\cdots$ is the unique positive root of octic equation 
$\chi_1(r) = 0$. 
Note that $a_8^{(1)} = 54.1247\cdots$ is nonzero at $r = r_1$. 
Thus Lemma \ref{lem:pde-octa} leads to Theorem \ref{thm:main-octa} 
for $k = 1$.  
\subsection{Face Problem} \label{subsec:face-octa}
For each index of $I_2$, the foot on the corresponding affine plane 
is given by  
\[
\rA^a\rB^b\rD^{abc} = (a q, \, b q, \, c p), \qquad
\rA^a\rC^c\rD^{abc} = (a q, \, b p, \, c q), \qquad
\rB^b\rC^c\rD^{abc} = (a p, \, b q, \, c q), 
\]
with $a, b, c = \pm$, where    
\[
p := \frac{r(3 - 2 r)}{3 (2 r^2 - 4 r + 3)}, \qquad  
q := \frac{r^2}{3 (2 r^2 - 4 r + 3)}. 
\]
\par
Observe that up to symmetries there are three types of 
vertex-edge-face flags: 
\begin{align*}
&(1) && \rA^+ \prec \rA^+\rB^+ \prec \rA^+\rB^+\rD^{+++} && 
(\VE 1) \,\, \& \,\, (\EF 1), \\[1mm] 
&(2) && \rA^+ \prec \rA^+\rD^{+++} \prec \rA^+\rB^+\rD^{+++} && 
(\VE 2) \,\, \& \,\, (\EF 2), \\[1mm] 
&(3) && \rD^{+++} \prec \rA^+\rD^{+++} \prec \rA^+\rB^+\rD^{+++} && 
(\VE 3) \,\, \& \,\, (\EF 2).  
\end{align*}
Formula \eqref{eqn:tau2-tetra} remains true for triakis octahedra if 
formulas in Table \ref{tab:tau2-tetra} are replaced by 
\begin{align*}
\tau_{1,m}^{(2)}(x) &= \sum h_m(\rA^a, \rA^a\rB^b, \rA^a\rB^b\rD^{abc}) + \cdots & &
\mbox{sum over flags of type (1)}, \\
\tau_{2,m}^{(2)}(x) &= \sum h_m(\rA^a, \rA^a\rD^{abc}, \rA^a\rB^b\rD^{abc}) + \cdots & & 
\mbox{sum over flags of type (2)}, \\
\tau_{3,m}^{(2)}(x) &= \sum h_m(\rD^{abc}, \rA^a\rD^{abc}, \rA^a\rB^b\rD^{abc}) + \cdots & &
\mbox{sum over flags of type (3)}.  
\end{align*}
\par
While the vertex-to-edge incidence numbers are given in 
\eqref{eqn:ve-octa}, the edge-to-face ones are  
\[
\ef_1 = \frac{3 - 2 r}{\sqrt{6 (2 r^2 - 4 r + 3)}}, \qquad 
\ef_2 = \frac{r(r - 1)}{\sqrt{(r^2 - 2 r + 3)(2 r^2 - 4 r + 3)}}.     
\]
So upon multiplying by a nonzero constant simultaneously, one may put 
\begin{align*}
\ve_1 \cdot \ef_1 &= \frac{3 - 2 r}{2 (2 r^2 - 4 r + 3)}, \qquad &   
\ve_2 \cdot \ef_2 &= \frac{(3 - r) r (r - 1)}{(r^2 - 2 r + 3)(2 r^2 - 4 r + 3)}, \\  
                  &                                   &
\ve_3 \cdot \ef_2 &= \frac{r^2(r - 1)^2}{(r^2 - 2 r + 3)(2 r^2 - 4 r + 3)}.   
\end{align*} 
\par
Putting all these informations together into formulas \eqref{eqn:coeff-octa}, 
one finds 
\begin{align*}
a_2^{(2)} &= \frac{8}{3}(r^2+2r+6) > 0, \qquad     
a_4^{(2)} = \frac{8}{27} \, \chi_2(r), \\ 
a_6^{(2)} &= \frac{8}{81}(16r^6+32r^5+24r^4-18r^3-54r^2+243) > 0, \\ 
a_8^{(2)} &= \frac{8}{2187}(16r^8+32r^7+24r^6-72r^5-324r^4-648r^3-486r^2+1458r+6561),  
\end{align*}
where $\chi_2(r)$ is the quartic polynomial defined in \eqref{eqn:quartic}. 
Observe that $a_2^{(2)}$ and $a_6^{(2)}$ are positive for every $r > 0$. 
Indeed the former is obvious and latter follows from the fact that 
$a_6^{(2)}$ $(r \ge 0)$ attains its minimum $22.0304\cdots > 0$ at 
$r = 0.743471\cdots$.    
Note that $a_4^{(2)} = 0$ if and only if $r$ is the unique positive root  
$r_2 = 1.82977\cdots$ of equation \eqref{eqn:quartic}, at which $a_8^{(2)} = 
13.2853\cdots$ is nonzero.  
Lemma \ref{lem:pde-octa} thus leads to Theorem \ref{thm:main-octa} for $k = 2$. 
\par 
By \cite[Theorem 2.2]{Iwasaki1} the volume problem ($k = 3$) has the same 
solution as the face problem ($k = 2$), since $P$ is isohedral. 
The proof of Theorem \ref{thm:main-octa} is now complete. 
\par
It is an interesting exercise left behind to deal with  a similar problem 
for isohedral triakis icosahedra.  
It would also be interesting if such special solids as discussed in this 
article appear in nature and physical sciences.  

\end{document}